\documentclass[10pt]{amsart}

\usepackage{amsmath,amssymb,amsthm}
\usepackage[margin=1in]{geometry}       
\usepackage{graphicx}                   
\usepackage{enumerate}                  
\usepackage{cite}                       

\usepackage{bookmark}                   

\usepackage{tikz}                       
\usetikzlibrary{decorations.pathreplacing}
\usetikzlibrary{patterns}
\tikzset{>=latex}

\usepackage{subcaption}                 
\usepackage{ctable}                     
\usepackage{multirow}                   
\usepackage{wrapfig}                    

\usepackage[T1]{fontenc}
\usepackage{charter}
\usepackage[expert]{mathdesign}

\usepackage[activate={true,nocompatibility},final,tracking=true,
kerning=true,spacing=true,factor=1100,stretch=10,shrink=10]{microtype}

\SetProtrusion{encoding={*},family={bch},series={*},size={6,7}}
               {1={ ,750},2={ ,500},3={ ,500},4={ ,500},5={ ,500},
                6={ ,500},7={ ,600},8={ ,500},9={ ,500},0={ ,500}}

\SetExtraKerning[unit=space]
    {encoding={*}, family={bch}, series={*}, size={footnotesize,small,normalsize}}
    {\textendash={400,400}, 
     "28={ ,150}, 
     "29={150, }, 
     \textquotedblleft={ ,150}, 
     \textquotedblright={150, }} 

\SetTracking{encoding={*}, family={bch}, series={*}, size={6,7,8,9,10,11,12}, shape=sc}{0}

\microtypecontext{spacing=nonfrench}

\newtheorem{theorem}{Theorem}
\newtheorem{lemma}[theorem]{Lemma}
\newtheorem{conjecture}[theorem]{Conjecture}
\newtheorem{question}[theorem]{Question}
\theoremstyle{definition}
\newtheorem{definition}[theorem]{Definition}
\theoremstyle{remark}
\newtheorem{remark}[theorem]{Remark}

\newcommand{\acts}{\curvearrowright}
\newcommand{\diam}{\mathrm{diam}}
\newcommand{\dist}{\mathrm{dist}}
\newcommand{\bndr}[1]{\partial{#1}}
\newcommand{\proj}{\mathrm{proj}}
\DeclareMathOperator{\inter}{\mathrm{int}}

\begin{document}

\title[Smooth orbit equivalence of multidimensional Borel flows]{Smooth orbit equivalence \\
  of multidimensional Borel flows}


\author{Konstantin Slutsky}
\address{Department of Mathematics\\
  Iowa State University\\
  411 Morrill Road \\
  Ames, IA 50011} \email{kslutsky@gmail.com}

\begin{abstract}
  Free Borel \( \mathbb{R}^{d} \)-flows are smoothly equivalent if there is a Borel bijection
  between the phase spaces that maps orbits onto orbits and is a \( C^{\infty} \)-smooth orientation
  preserving diffeomorphism between orbits.  We show that all free non-tame Borel
  \( \mathbb{R}^{d} \)-flows are smoothly equivalent in every dimension \( d \ge 2 \).  This answers
  a question of B.~Miller and C.~Rosendal.
\end{abstract}

\thanks{Konstantin Slutsky's research is partially supported by the ANR project AGRUME
  (ANR-17-CE40-0026).}

\maketitle

\section{Introduction}
\label{sec:introduction}

Let us begin by defining the notions mentioned in the title as well as the related concepts that are
needed to state the main results of our work. A \textbf{Borel flow} is a Borel action
\( \mathbb{R}^{d} \acts \Omega \) of the Euclidean group on a standard Borel space~\( \Omega \). An
action of \( \vec{r} \in \mathbb{R}^{d} \) upon \( x \in \Omega \) is denoted by \( x + \vec{r} \).
An \textbf{orbit equivalence} between two flows \( \mathbb{R}^{d} \acts \Omega \) and
\( \mathbb{R}^{d} \acts \Omega' \) is a Borel bijection \( \xi : \Omega \to \Omega' \) that sends
orbits onto orbits: \( \xi\bigl(x + \mathbb{R}^{d}\bigr) = \xi(x) + \mathbb{R}^{d} \) for all
\( x \in \Omega \); when such a map exists, we say that the flows are \textbf{orbit equivalent}. For
an action \( \mathbb{R}^{d} \acts \Omega \) we let \( E \) denote the corresponding orbit
equivalence relation: \( x E y \iff x + \mathbb{R}^{d} = y + \mathbb{R}^{d} \). When the action is
moreover free, \( \rho : E \to \mathbb{R}^{d} \) will stand for the associated \textbf{cocycle},
determined uniquely by the condition \( x + \rho(x, y) = y \). Given free Borel flows on phase
spaces \( \Omega \) and \( \Omega' \), any orbit equivalence \( \xi : \Omega \to \Omega' \) gives
rise to a map \( \alpha_{\xi} : \Omega \times \mathbb{R}^{d} \to \mathbb{R}^{d} \) defined by
\( \alpha_{\xi}(x, \vec{r}) = \rho\bigl(\xi(x), \xi(x + \vec{r})\bigr) \). A Borel orbit equivalence
\( \xi \) is said to be a \textbf{smooth equivalence} if
\( \alpha_{\xi}(x,\,\cdot\,) : \mathbb{R}^{d} \to \mathbb{R}^{d} \) is a \( C^{\infty} \)-smooth
orientation preserving diffeomorphism for all points \( x \in \Omega \).

\subsection{Prior work}
\label{sec:prior-work}

The concept of orbit equivalence originated in ergodic theory, where the set-up differs in two
essential aspects. First, one endows phase spaces of flows with probability measures. The flows are
then assumed to preserve (or to quasi-preserve) these measures. Likewise, orbit equivalence maps are
required to be at least quasi-measure-preserving. Second, all the properties of interest are
expected to hold up to a null set. For instance, an orbit equivalence map may mix elements between
orbits as long as this behavior is confined to a set of measure zero. The latter is a notable
relaxation of the Borel definition.

Smooth equivalence of one-dimensional flows, better known under the name of \textbf{time-change
  equivalence}, is closely connected to the notion of Kakutani equivalence of
automorphisms~\cite{MR14222}, and has been studied extensively since the pioneering works of
J.~Feldman~\cite{MR409763} and A.~Katok\cite{MR0412383,MR0442195}. An important milestone was the
monograph of D.~Ornstein, D.~Rudolph, and B.~Weiss~\cite{MR653094}, which showed, in particular,
that there is a continuum of pairwise time-change inequivalent ergodic measure-preserving flows. The
higher-dimensional case was considered by D.~Rudolph~\cite{MR536948}, where he found a striking
difference with the one-dimensional situation\textemdash{}all ergodic measure-preserving
\( \mathbb{R}^{d} \)-flows, \( d \ge 2 \), are smoothly equivalent. J.~Feldman obtained a similar
result for quasi-measure-preserving flows in~\cite{MR1113569,MR1163729}.

Note that in the definition of time-change equivalence, it is essential to allow for the orbit
equivalence maps to be quasi-measure-preserving even if all the flows are assumed to be
measure-preserving (see~\cite[Remark~4.5]{MR955378} regarding the connection between the
integrability of the cocycle as required in~\cite{MR0412383} and measure class preservation of the
orbit equivalence). This underlines the strength of Rudolph's result, as it is shown
in~\cite[Proposition~1.1]{MR536948} that in the dimensions \( d \ge 2 \) there is a single class of
ergodic measure-preserving flows under \emph{measure-preserving} smooth equivalence relation.
Further discussion on what restrictions on the orbit equivalence maps may produce finer equivalence
relations among higher-dimensional flows can also be found in~\cite{MR536948}.

In this paper, we are interested in the descriptive set-theoretical viewpoint. This means that
neither flows are assumed to preserve any measures (thus increasing the pool of flows to consider),
nor orbit equivalence maps have to be quasi-measure-preserving (which may potentially collapse
previously inequivalent flows into the same class). On the other hand, the necessity to run
constructions on all orbits may, in principle, increase the number of equivalence classes, as
complicated dynamics of a flow can be contained in a null set. All in all, this framework is in a
general position to the one of ergodic theory, and ahead of time, it is not clear how versatile
smooth equivalence will turn out to be. The key work that investigated the subject from this purely
Borel vantage point is the article by B.~Miller and C.~Rosendal~\cite{MR2578608}, where they studied
Kakutani equivalence of Borel automorphisms and classified one-dimensional flows up to descriptive
time-change equivalence.

\begin{theorem}[Miller--Rosendal~{\cite[Theorem B]{MR2578608}}]
  \label{thm:non-tame-flows-smooth-equivlaence}
  All non-tame\footnote{A flow \( \mathbb{R}^{d} \acts \Omega \) is \textbf{tame} if there is a
    Borel set \( S \subseteq \Omega \) that intersects every orbit of the flow in a single
    point. The term \textbf{smooth} is often used in the literature instead, but since we also work
    with diffeomorphisms, this word will be used in the traditional sense of differential geometry.
    Tame flows should be considered trivial in the context of the questions we are interested in
    this paper.}  free Borel \( \mathbb{R} \)-flows are smoothly equivalent.
\end{theorem}

As the one-dimensional case has been settled, they posed~\cite[Problem C]{MR2578608} the following
problem: ``Classify free Borel \( \mathbb{R}^{d} \)-actions on Polish spaces up to
(\( C^{\infty} \)-)time-change isomorphism.'' In other words, does the analog of D.~Rudolph and
J.~Feldman theorems hold? Are there two (non-tame) inequivalent free Borel
\( \mathbb{R}^{d} \)-flows for any \( d \ge 2 \)?

These and related topics were studied in~\cite{MR3984276}, where we showed that any two non-tame
free \( \mathbb{R}^{d} \)-flows, \( d \ge 2 \), are smoothly equivalent \emph{up to a compressible
  set}.  The method to prove this result was an expansion of the one used in~\cite{MR1113569}, and
such a statement is about as far as ergodic-theoretical methods can go, since a compressible set has
measure zero relative to any probability measure invariant under the flow.

\subsection{Main results}
\label{sec:main-results}

In the present work, we give a complete answer to Problem C of~\cite{MR2578608} by showing that all
non-tame free \( \mathbb{R}^{d} \)-flows, \( d \ge 2 \), are smoothly equivalent
(Theorem~\ref{thm:time-change-equivalence}).  Table~\ref{number-classes} provides a concise summary
and compares the number of classes up to smooth equivalence in ergodic theory and Borel dynamics.

\ctable[
label=number-classes,
pos=htb,
caption = {Number of classes of smooth orbit equivalence.},
mincapwidth =\textwidth,
]{ccc}{}{%
  \toprule
  & Ergodic Theory                & Borel Dynamics       \\
  \midrule
  \( d = 1 \)   & \( \mathfrak{c} \)-many~\cite{MR653094}          & one~\cite{MR2578608} \\
  \( d \ge 2 \) & one~\cite{MR536948,MR1113569} & one                  \\
  \bottomrule }

Many results in ergodic theory and Borel dynamics of \( \mathbb{Z}^{d} \) and \( \mathbb{R}^{d} \)
actions are based on the fact that such actions are (essentially) hyperfinite.  In ergodic theory,
this is manifested by a group of related theorems that usually go under the name of ``Rokhlin
Lemma''.  The key idea here is that one can find a measurable set that intersects every orbit of the
flow in a set of pairwise disjoint rectangles (more precisely, \( d \)-dimensional parallelepipeds).
Moreover, one often takes a sequence of such sets, where rectangles cohere and eventually cover all
the orbits (at least, up to a null set).  The details of the assumptions on such regions vary, but a
construction of this form is present in many arguments, including the references above.  The direct
analog of such a tower of coherent rectangular regions is not possible in Borel dynamics.  One,
therefore, has to rely on more complicated geometric shapes (see, for instance,~\cite[Theorem
1.16]{MR1900547} and~\cite{MR3359054}).

Our argument also requires regions witnessing hyperfiniteness.  The key property we need is for them
to be smooth disks.  S.~Gao, S.~Jackson, E.~Krohne, and B.~Seward~\cite{GaoJacksonDiskShapedRegions}
have shown the possibility to construct such regions for low-dimensional flows.  Their argument is
an elaboration of the orthogonal marker regions technique developed in~\cite{MR3359054}.  We take a
different path and build upon the approach presented by A.~Marks and S.~Unger
in~\cite[Appendix~A]{MR3702673}.  Section~\ref{sec:toasts} is devoted to these topics and it leads
to Theorem~\ref{thm:ball-shaped-toast} that shows the existence of such disk-shaped regions in all
dimensions.

In order to prove that all non-tame free \( \mathbb{R}^{d} \)-flows are smoothly equivalent, we
leverage the work of Miller--Rosendal that handles the case of \( d = 1 \).  To this end
Section~\ref{sec:equiv-to-special-flows} introduces the concept of a special flow, which is a type
of an \( \mathbb{R}^{d} \)-flow that is build over a one-dimensional flow in a very primitive way.
We show in Theorem~\ref{thm:special-flow-equivalence} that all \( \mathbb{R}^{d} \)-flows are
smoothly equivalent to a special flow.  This piece is the technical core of this paper.

Finally, in Section~\ref{sec:smooth-equiv-flows} we prove the main result on smooth equivalence of
free \( \mathbb{R}^{d} \)-flows (Theorem~\ref{thm:time-change-equivalence}) and conclude with some
remarks on its potential strengthening.

\subsection{Notations}
The following notations are used throughout the paper: \( B(R) \subseteq \mathbb{R}^{d} \) denotes a
ball of radius \( R \) centered at the origin; \( ||\,\cdot\,|| \) stands for the \( \ell^{2} \)-
norm in \( \mathbb{R}^{d} \); and \( \dist(\vec{r}, \vec{r}\,') \) refers to the Euclidean distance
in \( \mathbb{R}^{d} \). By a diffeomorphism we always mean a \( C^{\infty} \)-smooth orientation
preserving diffeomorphism. A smooth disk therefore refers to any compact region in
\( \mathbb{R}^{d} \) that is diffeomorphic to a ball. Interior of a set
\( F \subseteq \mathbb{R}^{d} \) is denoted by \( \inter F \), and \( \bndr F \) stands for the
boundary of \( F \). Given a Cartesian product \( X_{1} \times X_{2} \times \cdots \times X_{m} \),
\( \proj_{k} : \prod_{i =1}^{m}X_{i} \to X_{k} \) denotes the projection onto the
\( k^{\textrm{th}} \) coordinate, and, more generally, \( \proj_{[k,l]} \) will denote the
projection onto \( X_{k} \times \cdots \times X_{l} \), for \( k \le l \).

\subsection{Acknowledgement}
\label{sec:acknowledgement}

The author expresses his appreciation to Todor Tsankov for numerous helpful discussions on the
topic of this paper.

\section{Disk-Shaped Coherent Regions}
\label{sec:toasts}

We begin by stating the following classical fact from differential topology (see, for
instance,~\cite[Proposition 2.6]{MR1113569}), which will be used throughout the paper to justify the
existence of diffeomorphisms moving disks in a prescribed fashion.

\begin{lemma}[Extension Lemma]
  \label{lem:moving-smooth-disks}
  Let \( F \) and \( F' \) be smooth disks in \( \mathbb{R}^{d} \), \( d \ge 2 \), each containing
  \( m \) smooth disks in its interior: \( D_{1}, \ldots, D_{m} \subset \inter F \) and
  \( D_{1}', \ldots, D_{m}' \subset \inter F' \).  Suppose that disks \( D_{i} \) are pairwise
  disjoint and so are the disks \( D_{i}' \).  Any collection of orientation preserving
  diffeomorphisms \( \phi_{i} : D_{i} \to D_{i}' \) can be extended to an orientation preserving
  diffeomorphism \( \psi : F \to F' \).
\end{lemma}

Theorems in Borel dynamics of \( \mathbb{R}^{d} \) and \( \mathbb{Z}^{d} \) actions often rely on
variants of the hyperfiniteness construction.  Our argument is no exception, and this section gives
the specific version to be used later in Section~\ref{sec:equiv-to-special-flows}.  The cases of
\( d = 2 \) and \( d = 3 \) of Theorem~\ref{thm:ball-shaped-toast} are due to S.~Gao, S.~Jackson,
E.~Krohne, and B.~Seward; they are announced to appear in~\cite{GaoJacksonDiskShapedRegions}.  We
borrow the structure of our argument from A.~Marks and S.~Unger~\cite[Appendix~A]{MR3702673} and
supplement it with Lemma~\ref{lem:disk-separation} to get the desired shape of the regions for all
dimensions \( d \ge 2 \).

\begin{lemma}[Separation Lemma]
  \label{lem:disk-separation}
  Let \( 0 < R_{1} < R_{2} \) be positive reals and let
  \( D_{1}, \ldots, D_{n} \subset \mathbb{R}^{d} \), \( d \ge 2 \), be smooth pairwise disjoint
  disks of diameter \( \diam(D_{i}) < (R_{2}-R_{1})/2 \).  There exists a smooth disk \( F \) wedged
  between the two balls, \( B(R_{1}) \subseteq F \subseteq B(R_{2}) \), such that for each \( i \)
  either \( D_{i} \subseteq \inter F \) or \( F \cap D_{i} = \varnothing \).
\end{lemma}

\begin{wrapfigure}{r}{0.35\textwidth}
  \begin{center}
    \begin{tikzpicture}
      \draw[pattern=north east lines,pattern color=lightgray,rounded corners=1.4mm] (1.1,1.3) 
      -- (1.4, 0.6) -- (1.4, 0.1) -- (1, -0.2)
      -- (0.9, -0.5) -- (0.8, -1.0) -- (0.6, -1.4) -- (0.2, -1.5) -- (-0.1, -1.3)
      -- (-0.2, -1.0) -- (-0.8, -1.1) -- (-1.1, -0.8) -- (-1.3, -0.3) -- (-1.1, 0)
      -- (-0.9, 0.2) -- (-0.8, 0.5) -- (-0.8, 0.7) -- (-0.8, 1.0) -- (-0.9, 1.3)
      -- (-0.6, 1.6) -- (-0.4, 1.7) -- (-0.05, 1.8) -- (0.1, 1.4) -- (0.1, 0.9) -- (0.4, 0.7)
      -- (0.7, 1.1) -- cycle;
      \draw[thick] (0,0.05) circle (0.55cm);
      \draw[thick] (0,0.05) circle (1.85cm);
      \draw (0.2,0.3) node {\( R_{1} \)};
      \draw (1.6,1.4) node {\( R_{2} \)};
      \draw (1.0,-0.9) node {\( F \)};
      \draw (0,0) node {\includegraphics[width=4cm]{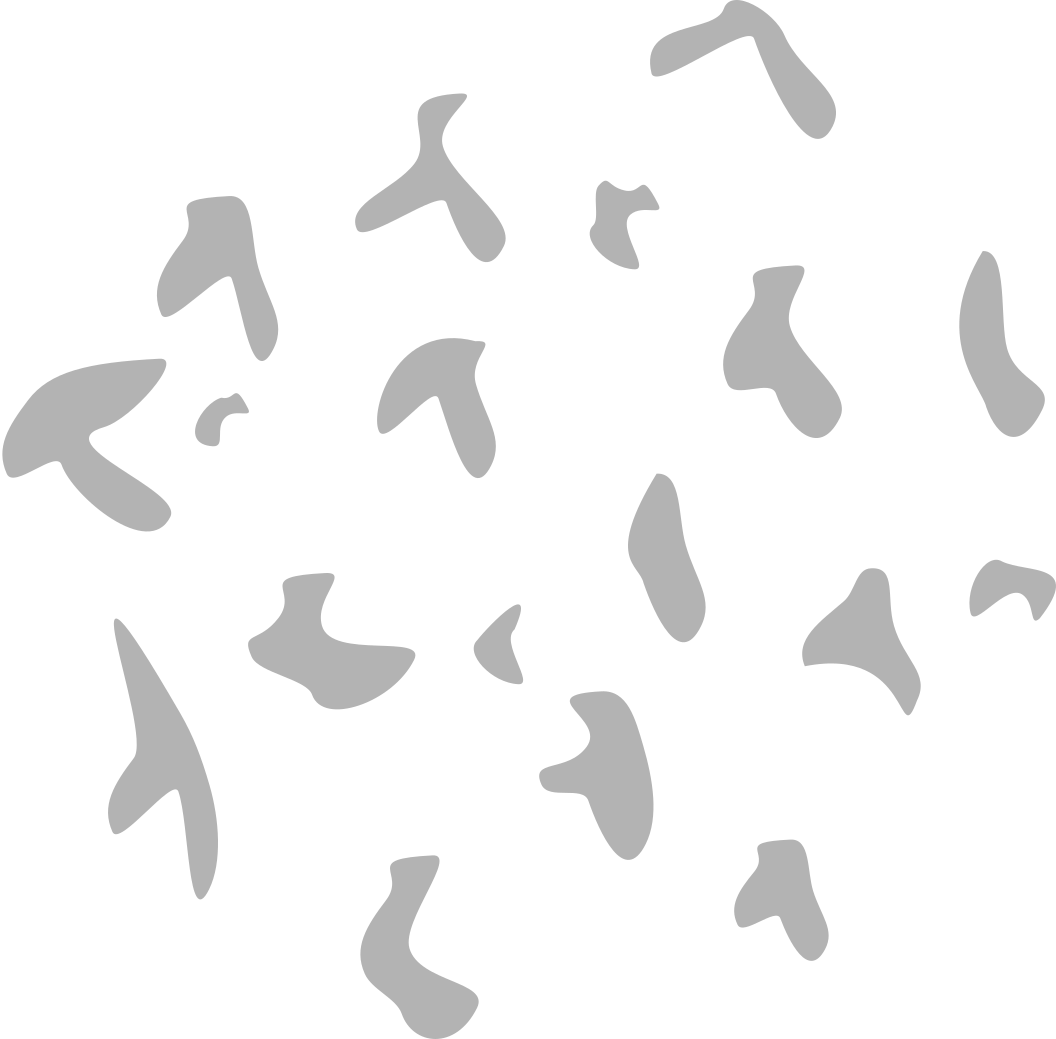}};
    \end{tikzpicture}
  \end{center}
  \caption{Separation Lemma.}
  \label{fig:disk-separation}
  \vspace{-0.1cm}
\end{wrapfigure}
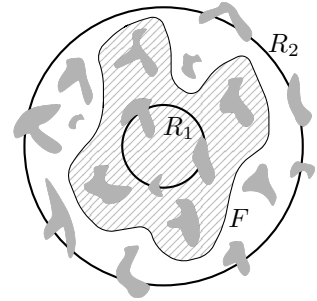

Figure~\ref{fig:disk-separation} illustrates the statement.  Disks \( D_{i} \) are marked in
gray and the required disk \( F \) is dashed.

\begin{proof}
  The proof is by induction on the number of disks \( n \).  The base case \( n = 0 \) is trivial, we
  argue the step from \( n-1 \) to \( n \).  If none of the disks \( D_{i} \) lie inside the open
  annulus \( A = \inter B(R_{2}) \setminus B(R_{1}) \), then the ball \( F = B((R_{2}+R_{1})/2) \)
  works.  Otherwise, select a ball \( D_{i_{0}} \subset A \).  By the inductive assumption there is a disk
  \( F' \) that fulfills the conclusions of the lemma for all disks \( D_{i} \), \( i \ne i_{0} \).
  We are done if also \( D_{i_{0}} \cap F' = \varnothing \) or \( D_{i_{0}} \subset \inter F' \), so
  assume otherwise (Figure~\ref{fig:step-one}).

  Find a smooth disk \( G \subset A \) that contains \( D_{i_{0}} \subset \inter G \) in its
  interior and does not intersect any other disk \( D_{i} \).  Pick a disk \( Z \subset \inter G \)
  that is disjoint from \( \bndr F' \) (Figure~\ref{fig:step-two}).  Such a disk can be found, since
  the boundary \( \bndr F' \) is nowhere dense.  Choose a diffeomorphism \( \psi \) supported on
  \( G \) such that \( \psi(D_{i_{0}}) = Z \).  Lemma~\ref{lem:moving-smooth-disks} may be used to
  justify the existence of such a diffeomorphism.  We have either
  \( \psi(D_{i_{0}}) \subset \inter F' \) or \( \psi(D_{i_{0}}) \cap F' = \varnothing \).  Set
  \( F = \psi^{-1}(F') \) (Figure~\ref{fig:step-three}).

  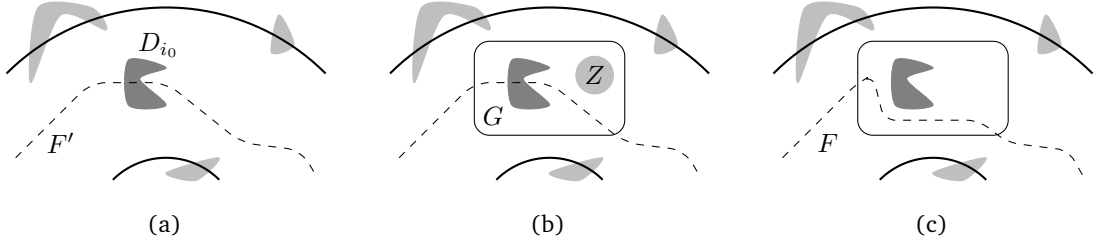
\begin{figure}[htb]
    \centering
    \begin{subfigure}{5cm}
      \begin{tikzpicture}
        \draw[white] (-2.5,0.5) rectangle (2.5,3.2);
        \filldraw[gray] plot [smooth cycle] coordinates {(-0.5, 1.7) (-0.5, 2.3) (0, 2.2) (-0.35,
          2.05) (0, 1.7)};
        \draw (-0.1, 2.5) node {\( D_{i_{0}} \)};
        \draw (-1.4, 1.2) node {\( F' \)};
        \draw[rounded corners=2mm, dashed] (-2,1) -- (-1,2) -- (0,2) -- (1,1.2) --
        (1.8, 1.1) -- (2, 0.75);
        \filldraw [lightgray] plot [smooth cycle] coordinates {(-1.8,2) (-1.7,3) (-1,3) (-0.8, 2.7)
        (-1.4, 2.9)};
        \filldraw [lightgray] plot [smooth cycle] coordinates {(0,0.8) (0.7, 1.0) (0.4, 0.7)};
        \filldraw [lightgray] plot [smooth cycle] coordinates {(1.5,2.9) (1.4, 2.4) (1.7, 2.5)};
        \draw [thick,  domain=45:135, samples=40] plot ({1*cos(\x)}, {1*sin(\x)});
        \draw [thick,  domain=45:135, samples=40] plot ({3*cos(\x)}, {3*sin(\x)});
      \end{tikzpicture}
      \caption{}
      \label{fig:step-one}
    \end{subfigure}
    \begin{subfigure}{5cm}
      \begin{tikzpicture}
        \draw[white] (-2.5,0.5) rectangle (2.5,3.2);
        \filldraw[gray] plot [smooth cycle] coordinates {(-0.5, 1.7) (-0.5, 2.3) (0, 2.2) (-0.35,
          2.05) (0, 1.7)};
        \draw[rounded corners=2mm] (-1, 1.3) rectangle (1,2.55);
        \draw (-0.75, 1.55) node {\( G \)};
        \filldraw[lightgray] (0.6, 2.1) circle (2.5mm);
        \draw (0.6, 2.1) node {\( Z \)};
        \draw[rounded corners=2mm, dashed] (-2,1) -- (-1,2) -- (0,2) -- (1,1.2) --
        (1.8, 1.1) -- (2, 0.75);
        \filldraw [lightgray] plot [smooth cycle] coordinates {(-1.8,2) (-1.7,3) (-1,3) (-0.8, 2.7)
        (-1.4, 2.9)};
        \filldraw [lightgray] plot [smooth cycle] coordinates {(0,0.8) (0.7, 1.0) (0.4, 0.7)};
        \filldraw [lightgray] plot [smooth cycle] coordinates {(1.5,2.9) (1.4, 2.4) (1.7, 2.5)};
        \draw [thick,  domain=45:135, samples=40] plot ({1*cos(\x)}, {1*sin(\x)});
        \draw [thick,  domain=45:135, samples=40] plot ({3*cos(\x)}, {3*sin(\x)});
      \end{tikzpicture}
      \caption{}
      \label{fig:step-two}
    \end{subfigure}
    \begin{subfigure}{5cm}
      \begin{tikzpicture}
        \draw[white] (-2.5,0.5) rectangle (2.5,3.2);
        \filldraw[gray] plot [smooth cycle] coordinates {(-0.5, 1.7) (-0.5, 2.3) (0, 2.2) (-0.35,
          2.05) (0, 1.7)};
        \draw[rounded corners=2mm] (-1, 1.3) rectangle (1,2.55);
        \draw[rounded corners=1.5mm, dashed] (-2,1) -- (-1,2) -- (-0.8, 2.1) -- (-0.65,1.5) -- (0.6, 1.5) --
        (1,1.2) -- (1.8, 1.1) -- (2, 0.75);
        \filldraw [lightgray] plot [smooth cycle] coordinates {(-1.8,2) (-1.7,3) (-1,3) (-0.8, 2.7)
        (-1.4, 2.9)};
        \filldraw [lightgray] plot [smooth cycle] coordinates {(0,0.8) (0.7, 1.0) (0.4, 0.7)};
        \filldraw [lightgray] plot [smooth cycle] coordinates {(1.5,2.9) (1.4, 2.4) (1.7, 2.5)};
        \draw (-1.4, 1.2) node {\( F \)};
        \draw [thick,  domain=45:135, samples=40] plot ({1*cos(\x)}, {1*sin(\x)});
        \draw [thick,  domain=45:135, samples=40] plot ({3*cos(\x)}, {3*sin(\x)});
      \end{tikzpicture}
      \caption{}
      \label{fig:step-three}
    \end{subfigure}
    \caption{Construction of a disk \( F \) that separates disks \( D_{i} \).}
    \label{fig:disk-separation-construction}
  \end{figure}

  Since \( \psi \) is supported on \( G \),
  both conditions \( F \cap D_{i} = \varnothing \) and \( D_{i} \subset \inter F \),
  \( i \ne i_{0} \), continue to hold whenever they did so for \( F' \) instead of \( F \).  By
  construction we now also have either \( F \cap D_{i_{0}} = \varnothing \) or
  \( D_{i_{0}} \subset \inter F \).
\end{proof}

Let \( \mathbb{R}^{d} \acts \Omega \) be a free Borel flow, and let \( E \) be its orbit equivalence
relation.  A set \( \mathcal{C} \subset \Omega \) is said to be
\begin{itemize}
\item \( R \)-\textbf{discrete}, where \( R \) is a positive real, if
  \( (c + B(R)) \cap (c' + B(R)) = \varnothing \) for all distinct \( c, c' \in \mathcal{C} \);
\item \textbf{discrete} if it is \( R \)-discrete for some \( R > 0 \);
\item \textbf{cocompact} if there exists \( R > 0 \) such that \( \mathcal{C} + B(R) = \Omega \);
\item \textbf{complete} if it intersects every orbit of the action;
\item a \textbf{cross section} if it is discrete and complete;
\item \textbf{on a rational grid} (or simply \textbf{rational}, for short) if
  \( \rho(c, c') \in \mathbb{Q}^{d} \) for all \( c, c' \in \mathcal{C} \) such that \( c E c' \).
\end{itemize}

\noindent We make use of the following result due to C.~M.~Boykin and S.~Jackson.

\begin{lemma}[Boykin--Jackson~\cite{MR2322367}, cf.~Lemma~A.2 of~\cite{MR3702673}]
  \label{lem:random-cross-sections}
  Let \( a_{1} < a_{2} < \cdots \) be an increasing sequence of natural numbers.  For any free Borel
  flow \( \mathbb{R}^{d} \acts \Omega \) there exists a sequence of \( a_{i} \)-discrete cocompact
  cross sections \( \mathcal{C}_{i} \) such that \( \bigcup_{i} \mathcal{C}_{i} \) is rational and
  for all \( \epsilon > 0 \), for every \( x \in \Omega \), there are infinitely many \( i \) such
  that \( ||\rho(x, c)|| < \epsilon a_{i} \) for some \( c \in \mathcal{C}_{i} \).
\end{lemma}

\begin{proof}
  The direct adaptation to \( \mathbb{R}^{d} \)-flows of the argument~\cite[Lemma~A.2]{MR3702673}
  (presented therein for~\( \mathbb{Z}^{d} \) actions) produces cross sections
  \( \mathcal{C}'_{i} \) that satisfy all the conclusions except possibly for
  \( \bigcup_{i} \mathcal{C}'_{i} \) being rational. As shown in~\cite[Lemma~2.3]{MR3984276}, there
  is a \emph{rational grid} for the flow, i.e., there is a complete rational set
  \( Q \subset \Omega \) invariant under the action of \( \mathbb{Q}^{d} \). Using Luzin-Novikov
  Theorem (see~\cite[18.14]{kechris_classical_1995}), one can find cross sections
  \( \mathcal{C}_{i} \subset Q \) and Borel bijections
  \( \zeta_{i} : \mathcal{C}_{i}' \to \mathcal{C}_{i} \) such that \( \bigcup_{i} \mathcal{C}_{i} \)
  is rational and for all \( c \in \mathcal{C}_{i}' \) one has \( c E\zeta_{i}(c) \) and
  \( ||\rho(c, \zeta_{i}(c))|| < 1 \). In other words, every element in \( \mathcal{C}_{i}' \) can
  be shifted by distance \( < 1 \) to ensure that all the cross sections are on the same rational
  grid. This argument is the content of~\cite[Lemma 2.4]{MR3984276}.

  The cross sections \( \mathcal{C}_{i} \) continue to be cocompact and still satisfy the key
  property that for every \( x \in \Omega \) and \( \epsilon > 0 \) there are infinitely many
  \( i \) with \( ||\rho(c,x)|| < \epsilon a_{i} \) for some \( c \in \mathcal{C}_{i} \). The only
  minor issue is that this modification reduces the discreteness parameter by~\( 1 \). Therefore if
  the original cross sections \( \mathcal{C}_{i}' \) were chosen to be \( (a_{i}+1) \)-discrete,
  then each of \( \mathcal{C}_{i} \) is guaranteed to be \( a_{i} \)-discrete.
\end{proof}

To formulate the next theorem we need an extra bit of notation.  Let
\( \mathbb{R}^{d} \acts \Omega \) be a free Borel flow.  For a set
\( W \subseteq \Omega \times \Omega \) and \( c \in \Omega \) we let \( W(c) \) denote the slice
over \( c \), i.e., \( W(c) = \{ x \in \Omega: (c,x) \in W \} \).  We also denote by
\( \widetilde{W} \) the set
\( \bigl\{(c, \vec{r}\,) \in \Omega \times \mathbb{R}^{d} : c+\vec{r} \in W(c) \bigr\} \).  Note that
\( \widetilde{W}(c) \) is the region of \( \mathbb{R}^{d} \) described by \( W(c) \), when \( c \)
is taken to be the origin of the coordinate system.

\begin{theorem}
  \label{thm:ball-shaped-toast}
  Let \( \mathbb{R}^{d} \acts \Omega \) be a free Borel flow and let \( E \) denote its orbit
  equivalence relation.  There exist cross sections~\( \mathcal{C}_{n} \) and Borel sets
  \( W_{n} \subseteq (\mathcal{C}_{n} \times \Omega) \cap E \) such that
  \( \bigcup_{n} \mathcal{C}_{n} \) is rational and for all \( n \in \mathbb{N} \):

  \begin{enumerate}[(i)]
  \item\label{item:regions-are-disks} \( \widetilde{W}_{n}(c) \) is a smooth disk for every
    \( c \in \mathcal{C}_{n} \).
  \item\label{item:regions-are-disjoint} Sets \( W_{n}(c) \), \( c \in \mathcal{C}_{n} \), are
    pairwise disjoint.
  \item\label{item:regions-are-coherent} For every \( c' \in \mathcal{C}_{m} \), \( m < n \), and
    every \( c \in \mathcal{C}_{n} \), either \( W_{m}(c') \cap W_{n}(c) = \varnothing \) or
    \(W_{m}(c') \subseteq W_{n}(c) \).  Moreover, in the latter case
    \( \rho(c,c') + \widetilde{W}_{m}(c') \) is contained in the interior of
    \( \widetilde{W}_{n}(c) \).
  \item\label{item:regions-exhaust-space} For all \( x \in \Omega \) and all compact
    \( K \subset \mathbb{R}^{d} \) there are \( m \) and \( c \in \mathcal{C}_{m} \) such that
    \( x + K \subseteq W_{m}(c) \).
  \item\label{item:regions-have-countably-many-shapes} There are smooth disks
    \( A_{n,k} \subseteq \mathbb{R}^{d} \), \( k \in \mathbb{N} \), and a Borel partition
    \( \mathcal{C}_{n} = \bigsqcup_{k \in \mathbb{N}} \mathcal{C}_{n,k} \) such that
    \[ W_{n} = \bigsqcup_{k} \bigl\{(c, c + \vec{r}\,) : c \in \mathcal{C}_{n,k}, \vec{r} \in
      A_{n,k}\bigr\} \quad \textrm{and} \quad \widetilde{W}_{n} = \bigsqcup_{k} \bigl\{(c, \vec{r}\,)
      : c \in \mathcal{C}_{n,k}, \vec{r} \in A_{n,k}\bigr\}.\]
  \end{enumerate}
\end{theorem}

\begin{proof} 
  Set \( a_{n} = 5^{n} \), and let \( \mathcal{C}_{n} \) be a sequence of cross sections produced by
  Lemma~\ref{lem:random-cross-sections}.  Note that \( \bigcup_{n} \mathcal{C}_{n} \) is guaranteed
  to be rational.  We construct a sequence of Borel sets
  \( W_{n} \subseteq (\mathcal{C}_{n} \times \Omega) \cap E \), which will also satisfy
  \begin{equation}
    \label{eq:1}
    B(a_{n}/2) \subseteq \widetilde{W}_{n}(c) \subseteq B(a_{n})  \quad \textrm{for all }
    c \in \mathcal{C}_{n}.
  \end{equation}
  This property will later be helpful in establishing item~\eqref{item:regions-exhaust-space}.
  
  For the base of the argument set
  \( W_{1} = \bigl\{ (c, c + \vec{r}\,) : c \in \mathcal{C}_{1}, \vec{r} \in B(a_{1}) \bigr\} \).
  Note that item~\eqref{item:regions-have-countably-many-shapes} holds with a trivial partition
  \( \mathcal{C}_{1,1} = \mathcal{C}_{1} \), \( \mathcal{C}_{1,j} = \varnothing \) for
  \( j \ge 2 \), and \( A_{1,1} = B(a_{1}) \).  Suppose now that \( W_{i} \) have been constructed
  for \( i < n \) and satisfy all the items of the theorem.  Cross section \( \mathcal{C}_{n} \) is
  \( a_{n} \)-discrete, so regions \( c + B(a_{n}) \) are pairwise disjoint as \( c \) ranges over
  \( \mathcal{C}_{n} \).

  For a given \( c \in \mathcal{\mathcal{C}}_{n} \) we consider regions \( W_{i}(c') \),
  \( i < n \), that intersect \( c + B(a_{n}) \) and that are not contained in a bigger such region.
  More formally, begin by choosing all the elements
  \( c_{1}^{n-1}, \ldots, c_{l_{n-1}}^{n-1} \in \mathcal{C}_{n-1} \) such that
  \( W_{n-1}(c_{j}^{n-1}) \cap (c + B(a_{n})) \ne \varnothing \); next, pick all
  \( c_{1}^{n-2}, \ldots, c_{l_{n-2}}^{n-2} \in \mathcal{C}_{n-2} \) such that
  \( W_{n-2}(c_{j}^{n-2}) \cap (c + B(a_{n})) \ne \varnothing \) and
  \( W_{n-2}(c_{j}^{n-2}) \cap W_{n-1}(c_{i}^{n-1}) = \varnothing \) for all
  \( 1 \le i \le l_{n-1} \); continue in the same fashion, terminating in a collection
  \( c^{1}_{1}, \ldots, c^{1}_{l_{1}} \in \mathcal{C}_{1} \) such that
  \( W_{1}(c_{j}^{1}) \cap (c + B(a_{n})) \ne \varnothing \) and
  \( W_{1}(c_{j}^{1}) \cap W_{k}(c_{i}^{k}) = \varnothing \) for all \( 2 \le k < n \), and all
  \( 1 \le i \le l_{k} \).  Note that in view of Eq.~\eqref{eq:1}, there can only be finitely many
  points \( c_{i}^{k} \) at each step.  Let
  \( c_{1}, \ldots, c_{l} \in \bigcup_{i < n} \mathcal{C}_{i} \) be an enumeration of the elements
  \( c_{i}^{k} \), \( 1 \le k < n \), \( 1 \le i \le l_{k} \), and let for \( 1 \le j \le l \), the
  number \( i(j)\) be such that \( c_{j} \in \mathcal{C}_{i(j)} \).

  Sets \( W_{i(j)}(c_{j}) \) are pairwise disjoint, and we therefore find ourselves in the set up of
  Lemma~\ref{lem:disk-separation}, where the ball \( B(a_{n}) \) interacts with a number of pairwise
  disjoint smooth disks \( \rho(c,c_{j}) + \widetilde{W}_{i(j)}(c_{j}) \), each having diameter
  \( \le 2*a_{n-1} < a_{n}/2 \).  Lemma~\ref{lem:disk-separation} claims that we can find a smooth
  disk \( F \) squeezed according to \( B(a_{n}/2) \subseteq F \subseteq B(a_{n}) \), and such that
  every region \( \rho(c,c_{j}) + \widetilde{W}_{i(j)}(c_{j}) \) is either contained in the interior
  of \( F \) or is disjoint from it.  Set \( W_{n}(c) = \{ c + \vec{r} : \vec{r} \in F \} \) and
  note that \( \widetilde{W}_{n}(c) \) fulfills Eq.~\eqref{eq:1}.

  We claim that this construction can be done in such a way that only countably many distinct shapes
  for \( F \) are used.  Indeed, the input to Lemma~\ref{lem:disk-separation}, which produced
  \( F \), is determined by the number \( l \) of regions \( W_{i(j)}(c_{j}) \) intersecting
  \( c + B(a_{n}) \), by the shape of these regions, and by their location relative to \( c \).
  Since the union \( \bigcup_{k} \mathcal{C}_{k} \) is rational, the vector
  \( (\rho(c, c_{1}), \ldots, \rho(c,c_{l})) \) is in \( \mathbb{Q}^{l} \).  By inductive
  assumption, for each \( c_{j} \in \mathcal{C}_{i(j)} = \bigsqcup_{k} \mathcal{C}_{i(j),k} \), there
  is some \( k(j) \in \mathbb{N} \) such that \( W_{i(j)}(c_{j}) = c_{j} + A_{i(j), k(j)} \) for a
  smooth disk \( A_{i(j), k(j)} \subseteq \mathbb{R}^{d} \).  Thus, the input to
  Lemma~\ref{lem:disk-separation} is uniquely determined by the tuple
  \[ \Bigl(l,\ \rho(c, c_{1}), \ldots, \rho(c, c_{l}),\ i(1), k(1), \ldots, i(l), k(l)\Bigr). \]
  There are only countably many such tuples and we can assume that the same disk \( F \) is used
  whenever the input tuple is the same. This guarantees compliance with
  item~\eqref{item:regions-have-countably-many-shapes}. Note also that such regions \( W_{n} \) are
  automatically Borel.
  
  It remains to verify that sets \( W_{n} \) satisfy the rest of the conclusions of the theorem.
  Item~\eqref{item:regions-are-disks} is fulfilled by the choice of \( F \).
  Item~\eqref{item:regions-are-disjoint} holds since \( \mathcal{C}_{n} \) is \( a_{n} \)-discrete
  and \( F \subseteq B(a_{n}) \).  Compliance with item~\eqref{item:regions-are-coherent} is the key
  property of the disk \( F \) produced by Lemma~\ref{lem:disk-separation}.

  We argue that item~\eqref{item:regions-exhaust-space} holds.  Pick a point \( x \in \Omega \) and
  a compact \( K \subset \mathbb{R}^{d} \).  Let \( n_{0} \) be so large that
  \( K \subseteq B(a_{n_{0}}/4) \).  According to the property of cross sections
  \( \mathcal{C}_{n} \) guaranteed by Lemma~\ref{lem:random-cross-sections}, for
  \( \epsilon = 1/4 \) there exists \( n_{1} \ge n_{0} \) such that \( ||\rho(x, c)|| < a_{n_{1}}/4 \)
  for some \( c \in \mathcal{C}_{n_{1}} \), i.e.,\ \( x \in c + B(a_{n_{1}}/4) \).  We
  therefore have
  \[ x + K \subseteq x + B(a_{n_{0}}/4) \subseteq c + B(a_{n_{1}}/4) + B(a_{n_{0}}/4) \subseteq c
    + B(a_{n_{1}}/2) \subseteq W_{n_{1}}(c), \]
  where the last inclusion follows from Eq.~\eqref{eq:1}.
\end{proof}

In the proof above we chose a family of pairwise disjoint regions \( W_{i(j)}(c_{j}) \) that
intersect \( c + B(a_{n}) \).  In the sequel, we will need a similar family of subregions of a
region \( W_{n}(c) \).  The following lemma and definition isolate the relevant notion.

\begin{lemma}
  \label{lem:maximal-sub-family}
  Let \( \mathcal{C}_{n} \) and \( W_{n} \), \( n \in \mathbb{N} \), be as in
  Theorem~\ref{thm:ball-shaped-toast}.  For each \( n \) and each \( c \in \mathcal{C}_{n} \) there
  exists a family \( c_{1}, \ldots, c_{l} \in \bigcup_{i < n} \mathcal{C}_{i} \) 
  such that for \( i(j) \) given by the condition \( c_{j} \in \mathcal{C}_{i(j)} \) one has
  \begin{enumerate}[(i)]
    \item\label{item:max-sub-family-containment} \( W_{i(j)}(c_{j}) \subseteq W_{n}(c) \) for all
      \( 1 \le j \le l \);
    \item\label{item:max-sub-family-disjoint} sets \( W_{i(j)}(c_{j}) \) are pairwise disjoint for
      \( 1 \le j \le l \);
    \item\label{item:max-sub-family-maximal} for any \( m < n \) and \( c' \in \mathcal{C}_{m} \)
      such that \( W_{m}(c') \subseteq W_{n}(c) \) there exists \( 1 \le j \le l \) such that
      \( W_{m}(c') \subseteq W_{i(j)}(c_{j}) \).
  \end{enumerate}
\end{lemma}

\begin{proof}
  Just like in the proof of Theorem~\ref{thm:ball-shaped-toast}, let
  \( c_{1}^{n-1}, \ldots, c_{l_{n-1}}^{n-1} \in \mathcal{C}_{n-1} \) be all the elements (if any)
  such that \( W_{n-1}(c_{j}^{n-1}) \subseteq W_{n}(c) \).  In view
  of~\ref{thm:ball-shaped-toast}\eqref{item:regions-are-disjoint}, sets \( W_{n-1}(c_{j}^{n-1}) \)
  are pairwise disjoint.  Pick all the elements
  \( c_{1}^{n-2}, \ldots, c_{l_{n-2}}^{n-2} \in \mathcal{C}_{n-2} \) satisfying
  \( W_{n-2}(c_{j}^{n-2}) \subseteq W_{n}(c) \), but \( W_{n-2}(c_{j}^{n-2}) \) is disjoint from all
  \( W_{n-1}(c_{i}^{n-1}) \), \( 1 \le i \le l_{n-1} \).  Note that
  by~\ref{thm:ball-shaped-toast}\eqref{item:regions-are-coherent} the latter is equivalent to
  saying that \( W_{n-2}(c_{j}^{n-2}) \) is not contained in any of \( W_{n-1}(c_{k}^{n-1}) \),
  \( 1 \le k \le l_{n-1} \).

  One continues in the same fashion.  At step \( k \) we pick elements
  \( c_{1}^{n-k}, \ldots, c_{l_{n-k}}^{n-k} \in \mathcal{C}_{n-k}\) that are contained in
  \( W_{n}(c) \) and are disjoint from all the sets \( W_{n-j}(c_{i}^{n-j}) \), \( 1 \le j < k \),
  \( 1 \le i \le l_{n-j} \), constructed at the previous steps.  The process terminates with the
  selection of elements \( c_{1}^{1}, \ldots, c_{l_{1}}^{1} \in \mathcal{C}_{1} \).
  
  The points \( c_{j}^{k} \), \( 1 \le k < n \), \( 1 \le j \le l_{k} \), satisfy the conditions of
  this lemma.  Items~\eqref{item:max-sub-family-containment}
  and~\eqref{item:max-sub-family-disjoint} are evident, and~\eqref{item:max-sub-family-maximal}
  follows from the observation that if \( W_{m}(c') \subseteq W_{n}(c) \) was not picked during the
  construction, then it had to intersect some set \( W_{k}(c_{j}^{k}) \) for an element
  \( c_{j}^{k} \), \( k > m \), picked earlier.  By the
  condition~\ref{thm:ball-shaped-toast}\eqref{item:regions-are-coherent} this means
  \( W_{m}(c') \subseteq W_{k}(c_{j}^{k}) \) as desired.
\end{proof}

\begin{definition}
  \label{def:maximal-sub-family}
  A family of regions \( W_{i(j)}(c_{j}) \) satisfying the conclusions of
  Lemma~\ref{lem:maximal-sub-family} is called \textbf{a maximal family of subregions}
  of~\( W_{n}(c) \).
\end{definition}

\begin{remark}
  \label{rem:max-subfamily-is-unique}
  It is easy to check that the maximal family of subregions of any \( W_{n}(c) \) is necessarily
  unique, but this will not play a role in our arguments.
\end{remark}

\begin{lemma}
  \label{lem:toast-properties}
  Let \( \mathcal{C}_{n} \) and \( W_{n} \), \( n \in \mathbb{N} \), be as in
  Theorem~\ref{thm:ball-shaped-toast}. For every \( m_{1} \in \mathbb{N} \) and
  \( c_{1} \in \mathcal{C}_{m_{1}} \) there exist a sequence of integers
  \( m_{1} < m_{2} < m_{3} < \cdots \) and elements \( c_{j} \in \mathcal{C}_{m_{j}} \) such that
  the regions \( W_{m_{j}}(c_{j}) \) satisfy the inclusions
  \( W_{m_{j}}(c_{j}) \subseteq W_{m_{j+1}}(c_{j+1}) \) for all \( 1 \le j < \infty \).
\end{lemma}

\begin{proof}
  The set \( \widetilde{W}_{m_{1}}(c_{1}) \) is a disk
  by~\ref{thm:ball-shaped-toast}\eqref{item:regions-are-disks}, and in particular it is a compact
  region in \( \mathbb{R}^{d} \). We may therefore pick a compact \( K \subseteq \mathbb{R}^{d} \)
  such that the inclusion \( \widetilde{W}_{m_{1}}(c_{1}) \subset K \) is proper.
  By~\ref{thm:ball-shaped-toast}\eqref{item:regions-exhaust-space} there exists some \( m_{2} \) and
  \( c_{2} \in \mathcal{C}_{m_{2}} \) such that
  \( W_{m_{1}}(c_{1}) \subset c_{1} + K \subseteq W_{m_{2}}(c_{2}) \).
  Items~\ref{thm:ball-shaped-toast}\eqref{item:regions-are-coherent}
  and~\ref{thm:ball-shaped-toast}\eqref{item:regions-are-disjoint} guarantee that
  \( m_{2} > m_{1} \). The same choice can now be iterated to construct the desired sequence
  \( m_{1} < m_{2} < m_{3} < \cdots \) and elements \( c_{j} \in \mathcal{C}_{m_{j}} \).
\end{proof}

\section{Equivalence to Special Flows}
\label{sec:equiv-to-special-flows}

One of the simplest ways to construct an \( \mathbb{R}^{d} \)-flow is to start with an
\( \mathbb{R} \)-flow on some standard Borel space~\( \Omega_{1} \) and define the action
\( \mathbb{R}^{d} \acts \Omega_{1} \times \mathbb{R}^{d-1} \) by
\[ \Omega_{1} \times \mathbb{R}^{d-1} \ni (y, \vec{q}\,) + (r_{1}, \ldots, r_{d}) = (y + r_{1},
  \vec{q} + (r_{2}, \ldots, r_{d})). \]
We say that a flow \( \mathbb{R}^{d} \acts \Omega \) is \textbf{special} if it is isomorphic to a
flow of the form above.  This is an ad hoc notion, which we use to reduce smooth equivalence of
multidimensional flows to the one dimensional situation.  Our goal in this section is to show that
every free Borel \( \mathbb{R}^{d} \)-flow is smoothly equivalent to a special one.  The argument
goes through a sequence of lemmas, and we begin by establishing some common notation.

Throughout the section we fix a free Borel \( \mathbb{R}^{d} \)-flow \( \mathfrak{F} \),
\( d \ge 2 \), let \( \mathcal{C}_{n} \) be the cross sections and
\( W_{n} \subseteq \mathcal{C}_{n} \times \Omega \) be the corresponding regions produced by
Theorem~\ref{thm:ball-shaped-toast}. Let \( V_{n} \subseteq \Omega \) denote the projection of
\( W_{n} \) onto the second coordinate, and let \( \pi_{n} : V_{n} \to \mathcal{C}_{n} \) be defined
by the condition \( (\pi_{n}(x), x) \in W_{n} \) for all \( x \in V_{n} \). Note that \( V_{n} \) is
Borel as \( \proj_{2} : W_{n} \to V_{n} \) is injective
by~\ref{thm:ball-shaped-toast}\eqref{item:regions-are-disjoint} and \( \pi_{n} \) is Borel since its
graph is the flip of \( W_{n} \). Define for \( m < n \) sets
\[ P_{m,n} = \bigl\{(c',c) \in \mathcal{C}_{m} \times \mathcal{C}_{n} : W_{m}(c') \subseteq W_{n}(c)
  \bigr\}, \]
which encode regions of the level \( m \) inside a given region of the level \( n \).  Sets
\( \widetilde{W}_{n}(c) \) are smooth disks, and our first lemma shows that specific diffeomorphisms
onto balls can be chosen to cohere across levels.

\begin{lemma}
  \label{lem:coherent-diffs}
  There exist radius maps \( t_{n} : \mathcal{C}_{n} \to \mathbb{R}^{> 0} \), ``diffeomorphism''
  functions \( \phi_{n} : V_{n} \to \mathbb{R}^{d} \), and shift maps
  \( s_{m,n} : P_{m,n} \to \mathbb{R}^{\ge 0} \) subject to the following conditions to be valid for
  all \( m < n \), and all \( (c', c) \in P_{m,n} \)\,:
  \begin{enumerate}[(i)]
    \item\label{item:radii-grow} \( t_{m}(c') \ge m \);
    \item\label{item:phi-is-a-diffeormophism}
      \( \widetilde{W}_{m}(c') \ni \vec{r} \mapsto \phi_{m}(c' + \vec{r}) \in B(t_{m}(c')) \) is
      a \( C^{\infty} \) orientation preserving diffeomorphism onto the ball \( B(t_{m}(c')) \);
    \item\label{item:coherence-of-steps} \( \phi_{n} (x) = \phi_{m}(x) + \vec{s}_{m,n}(c',c) \) for
      all \( x \in W_{m}(c') \), where
      \( \vec{s}_{m,n}(c',c) = s_{m,n}(c',c) \times \vec{0}^{d-1} \);
    \item\label{item:far-from-boundary} \( t_{m}(c') + s_{m,n}(c',c) \le t_{n}(c) -1 \);
    \item\label{item:countably-many-cases} there is a Borel partition
      \( \mathcal{C}_{m} = \bigsqcup_{k} \mathcal{C}_{m,k} \) such that
      \( \widetilde{W}_{m}(c_{1}) = \widetilde{W}_{m}(c_{2}) \), \( t_{m}(c_{1}) = t_{m}(c_{2}) \),
      and \( \phi_{m}(c_{1} + \vec{r}) = \phi_{m}(c_{2} + \vec{r}) \) for all
      \( c_{1},c_{2} \in \mathcal{C}_{m,k} \) and all \( \vec{r} \in \widetilde{W}_{m}(c_{1}) \).
  \end{enumerate}
\end{lemma}

The meaning of these conditions is as follows. Item~\eqref{item:radii-grow} ensures that radii go to
infinity as \( m \to \infty \). According to item~\eqref{item:phi-is-a-diffeormophism}, each map
\( \phi_{m} \) encodes a family of diffeomorphisms, one for each \( c' \in \mathcal{C}_{m} \).
Formally speaking, these diffeomorphisms are maps from \( \widetilde{W}_{m}(c') \) onto
\( B(t_{m}(c')) \). However, we will occasionally abuse the language by calling the map
\( \phi_{m}|_{W_{m}(c')} : W_{m}(c') \to B(t_{m}(c')) \) a diffeomorphism.

Item~\eqref{item:coherence-of-steps} postulates that \( \phi_{n} \) extends \( \phi_{m} \) up to a
translation of the range along the \( x \)-axis, where the translation value is constant over each
\( W_{m}(c') \) region and is equal to \( s_{m,n}(c',c) \). Condition~\eqref{item:far-from-boundary}
is a reformulation of the inequality
\[ \dist\bigl(\phi_{m}(W_{m}(c')) + \vec{s}_{m,n}(c',c),\ \bndr B(t_{n}(c))\bigr) \ge 1. \] It means
that disks \( \phi_{m}(W_{m}(c')) + \vec{s}_{m,n}(c',c) \) are at least one unit of distance away
from the boundary of \( B(t_{n}(c)) \) (see Figure~\ref{fig:alignment-of-images}). Similarly
to Theorem~\ref{thm:ball-shaped-toast}\eqref{item:regions-have-countably-many-shapes},
item~\eqref{item:countably-many-cases} says that we need to consider only countably many different
diffeomorphisms \( \phi_{m}|_{W_{m}(c')} \). The only purpose of this property is to make it easy
for us to argue that the flow \( \mathfrak{F}' \), which will be constructed later in this section,
is Borel.

\begin{figure}[hbt]
  \centering
  \begin{tikzpicture}
    \draw[rounded corners=4mm] (0,0) rectangle (8,5);
    \draw[rounded corners=2mm,fill=gray!20] (1.5,3) rectangle (3.5,4);
    \draw[rounded corners=4mm,fill=gray!20] (0.8,0.8) rectangle (3.3,2.3);
    \draw[rounded corners=10mm,fill=gray!20] (4.5,2) rectangle (7,4.5);
    \draw[rounded corners=1mm,fill=gray!20] (5.5,0.5) rectangle (6.5,1.3);

    \draw[rounded corners=3mm, dashed,very thick] (0,3.5) -- (1.5,3.5) --
    (3.5,3.5) -- (4.5,3.5) -- (5.5, 2.9) -- (7,3.7) -- (7.5,3.5) --
    (7.5, 1.7) -- (4, 1.7) -- (4, 2.6) -- (0.4, 2.6) -- (0.4, 1.2) --
    (2.2, 1.5) -- (4.3, 1.2) -- (4.9, 1.0) -- (6, 1.0) -- (8, 0.3);

    \draw[line width=0.7mm] (1.5,3.5) -- (3.5,3.5);

    \begin{scope}
      \clip[rounded corners=10mm](4.5,2) rectangle (7,4.5);
      \draw[rounded corners=3mm, line width=0.7mm] (0,3.5) -- (1.5,3.5) --
      (3.5,3.5) -- (4.5,3.5) -- (5.5, 2.9) -- (7,3.7) -- (7.5,3.5) --
      (7.5, 1.7) -- (4, 1.7) -- (4, 2.6) -- (0.4, 2.6) -- (0.4, 1.2) --
      (4.3, 1.2) -- (4.9, 1.0) -- (6, 1.0) -- (8, 0.3);
    \end{scope}
    \begin{scope}
      \clip[rounded corners=4mm] (0.8,0.8) rectangle (3.3,2.3);
      \draw[rounded corners=3mm, line width=0.7mm] (0,3.5) -- (1.5,3.5) --
      (3.5,3.5) -- (4.5,3.5) -- (5.5, 2.9) -- (7,3.7) -- (7.5,3.5) --
      (7.5, 1.7) -- (4, 1.7) -- (4, 2.6) -- (0.4, 2.6) -- (0.4, 1.2) --
      (2.2, 1.5) -- (4.3, 1.2) -- (4.9, 1.0) -- (6, 1.0) -- (8, 0.3);
    \end{scope}
    \begin{scope}
      \clip[rounded corners=1mm] (5.5,0.5) rectangle (6.5,1.3);
      \draw[rounded corners=3mm, line width=0.7mm] (0,3.5) -- (1.5,3.5) --
      (3.5,3.5) -- (4.5,3.5) -- (5.5, 2.9) -- (7,3.7) -- (7.5,3.5) --
      (7.5, 1.7) -- (4, 1.7) -- (4, 2.6) -- (0.4, 2.6) -- (0.4, 1.2) --
      (4.3, 1.2) -- (4.9, 1.0) -- (6, 1.0) -- (8, 0.3);
    \end{scope}
    \draw (5.75, 4) node {\( W_{i(j)}(c_{j}) \)};
    \draw (5.75, 2.7) node {\( L_{i(j)}(c_{j}) \)};
    \draw (4.3, 0.7) node {\( L_{n}(c) \)};
    \draw (0.8, 4.5) node {\( W_{n}(c) \)};
  \end{tikzpicture}
  \caption{Region \( W_{n}(c) \), containing four subregions \( W_{i(j)}(c_{j}) \) marked in
    gray.  Each of the subregions has a line segment \( L_{i(j)}(c_{j}) \), which are all
    contained inside \( L_{n}(c) \). }
  \label{fig:inductive-step}
\end{figure}
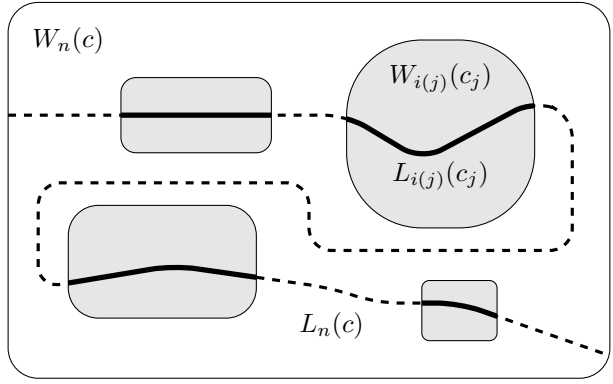

\begin{proof}[Proof of Lemma~\ref{lem:coherent-diffs}]
  The construction goes by induction on \( n \), and we begin with its base. Set \( t_{1}(c) = 1 \)
  for all \( c \in \mathcal{C}_{1} \). By item~\eqref{item:regions-are-disks} of
  Theorem~\ref{thm:ball-shaped-toast}, each region \( \widetilde{W}_{1}(c) \),
  \( c \in \mathcal{C}_{1} \), is a smooth disk. So for \( \phi_{1} : V_{1} \to B(1) \) we pick any
  map satisfying~\eqref{item:phi-is-a-diffeormophism} and~\eqref{item:countably-many-cases}, which
  can be done since there are only countably many shapes \( \widetilde{W}_{1}(c) \)
  by~\ref{thm:ball-shaped-toast}\eqref{item:regions-have-countably-many-shapes}.
  
  For the inductive step consider a region \( W_{n}(c) \). We pick points
  \( c_{1}, \ldots, c_{l} \in \bigcup_{i < n}\mathcal{C}_{i} \), and integers \( i(j) \) defined by
  \( c_{j} \in \mathcal{C}_{i(j)} \), that correspond to a maximal family of subregions of
  \( W_{n}(c) \) as per Lemma~\ref{lem:maximal-sub-family}. For such points \( c_{j} \) we have
  \( \rho(c, c_{j}) + \widetilde{W}_{i(j)}(c_{j}) \subset \inter\widetilde{W}_{n}(c) \), as
  guaranteed by~\ref{thm:ball-shaped-toast}\eqref{item:regions-are-coherent}. An example of such a
  region \( W_{n}(c) \) is shown in Figure~\ref{fig:inductive-step} on
  page~\pageref{fig:inductive-step}. By inductive assumption regions
  \( \widetilde{W}_{i(j)}(c_{j}) \) are diffeomorphic to balls \( B(t_{i(j)}(c_{j})) \) via the
  diffeomorphisms \( \vec{r} \mapsto \phi_{i(j)}(c_{j} + \vec{r}) \). We shift these balls along the
  \( x \)-axis to make them disjoint, and view them inside a sufficiently large ball in
  \( \mathbb{R}^{d} \) (see Figure~\ref{fig:alignment-of-images}). More specifically,
  let
  \[ W_{i(j),n}' = \{ (c',x) \in W_{i(j)} : c' \in \proj_{1}(P_{i(j),n}) \textrm{ and
    } c' \not \in \proj_{1}(P_{i(j), k}) \textrm{ for all } i(j) < k < n \}, \] and put
  \( V_{i(j),n}' = \proj_{2}(W'_{i(j),n}) \); in other words, \( W_{i(j), n}' \) consists of those
  regions \( W_{i(j)}(c') \) for which \( n \) is the smallest index to satisfy
  \( W_{i(j)}(c') \subset W_{n}(c) \) for some \( c \in \mathcal{C}_{n} \). Set for
  \( 1 \le j \le l \)
  \[ s_{i(j),n}(c_{j},c) = (j-1) + t_{i(j)}(c_{j}) + 2 \sum_{1 \le k < j} t_{i(k)}(c_{k}), \] and
  consider the map \( \phi'_{i(j)} : V_{i(j),n}' \to \mathbb{R}^{d} \) to be given for
  \( x \in W_{i(j)}(c_{j}) \) by
  \( \phi'_{i(j)}(x) = \phi_{i(j)}(x) + \vec{s}_{i(j), n}(c_{j},c) \). Note that restrictions
  \( \phi'_{i(j)}|_{W_{i(j)}(c_{j})} \) are diffeomorphisms onto disks
  \( B(t_{i(j)}(c_{j})) + \vec{s}_{i(j),n}(c_{j},c) \).
  The radius \( t_{n}(c) \) is taken to be sufficiently large to contain these disks:
  \( t_{n}(c) = \max\bigl\{n,\ s_{i(l), n}(c_{l},c) + t_{i(l)}(c_{l}) + 1 \bigr\} \). This ensures
  that images \( \phi'_{i(j)}(W_{i(j)}(c_{j})) \) are inside \( B(t_{n}(c)) \), and are furthermore
  at least \( 1 \) unit of distance away from its boundary, which yields
  item~\eqref{item:far-from-boundary}. Item~\eqref{item:radii-grow} also continues to be satisfied
  by this choice of \( t_{n}(c) \).

  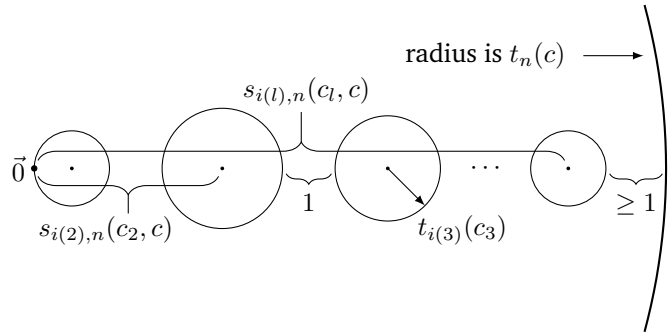
\begin{figure}[htb]
    \centering
    \begin{tikzpicture}
      \filldraw (0,0) circle (1pt) node[anchor=east] {\( \vec{0} \)};
      \draw (0.5,0) circle (0.5cm);
      \filldraw (0.5,0) circle (0.5pt);
      \draw (2.5,0) circle (0.8cm);
      \filldraw (2.5,0) circle (0.5pt);
      \draw (4.7,0) circle (0.7cm);
      \filldraw (4.7,0) circle (0.5pt);
      \draw (6,0) node {\( \cdots \)};
      \draw (7.1,0) circle (0.5cm);
      \filldraw (7.1,0) circle (0.5pt); \draw
      [thick,domain=-15:15] plot ({8.4*cos(\x)}, {8.4*sin(\x)}); \draw (6,1.5) node {radius is
        \( t_{n}(c) \) }; \draw[->] (7.3, 1.5) -- (8.1, 1.5);
      
      \draw [decorate,decoration={brace,mirror,amplitude=4pt},yshift=-2pt] (3.35,0) -- (3.95, 0)
      node [black, midway, anchor=north, yshift=-1.5mm] {\( 1 \)};
      
      \draw [decorate,decoration={brace,mirror,amplitude=7pt},yshift=-3pt] (0.05,0) -- (2.45, 0)
      node [black, midway, anchor=north, yshift=-4mm,xshift=-3mm] {\( s_{i(2),n}(c_{2},c) \)};
      \draw (1.25, -0.35) -- (1.25, -0.6);
      
      \draw [decorate,decoration={brace,amplitude=7pt},yshift=3pt] (0.05,0) -- (7.05, 0) node
      [black, midway, anchor=south, xshift=2pt, yshift=6mm] {\( s_{i(l),n}(c_{l},c) \)};
      \draw (3.55, 0.35) -- (3.55, 0.75);
      
      \draw[->] (4.7,0) -- ({4.7 + 0.7*cos(45)}, {-0.7*sin(45)}) node[anchor=north, xshift=5mm]
      {\( t_{i(3)}(c_{3}) \)};
      
      \draw [decorate,decoration={brace,mirror,amplitude=4pt},yshift=-2pt] (7.65,0) -- (8.35, 0)
      node [black, midway, anchor=north, yshift=-1.5mm] {\( \ge 1 \)};
    \end{tikzpicture}
    \caption{Alignment of disks
      \( \phi_{i(j)}(W_{i(j)}(c_{j})) +\vec{s}_{i(j),n}(c_{j},c) \).}
    \label{fig:alignment-of-images}
  \end{figure}
  
  Extension Lemma~\ref{lem:moving-smooth-disks} can now be applied to diffeomorphisms
  \( \rho(c, c_{j}) + \widetilde{W}_{i(j)}(c_{j}) \ni \vec{r} \mapsto \phi'_{i(j)}(c + \vec{r}) \),
  since they are defined on disjoint disks \( \rho(c,c_{j}) + \widetilde{W}_{i(j)}(c_{j}) \) and
  have disjoint images \( B(t_{i(j)}(c_{j})) + \vec{s}_{i(j), n}(c_{j},c) \), \(1 \le j \le l \).
  All the domains of these maps lie in the interior of the disk \( \widetilde{W}_{n}(c) \), while
  the images are subsets of \( \inter B(t_{n}(c)) \). We therefore can find a common extension to a
  diffeomorphism \( \phi_{n}|_{W_{n}(c)} : W_{n}(c) \to B(t_{n}(c)) \) that satisfies
  \[ \phi_{n}(x) = \phi_{i(j)}'(x) = \phi_{i(j)}(x) + \vec{s}_{i(j),n}(c_{j},c) \] for all
  \( x \in W_{i(j)}(c_{j}) \textrm{ and all } 1 \le j \le l \), thus
  implying~\eqref{item:coherence-of-steps}.

  We are not quite done yet though. The construction above defined diffeomorphisms \( \phi_{n} \)
  and radii \( t_{n}(c) \) for all \( c \in \mathcal{C}_{n} \), but shifts \( s_{m,n}(c',c) \) are
  currently defined only for those \( c' \in \mathcal{C}_{m} \) that belong to the maximal family of
  subregions of \( W_{n}(c) \). Nonetheless, values \( s_{m,n}(c',c) \) satisfying
  item~\eqref{item:coherence-of-steps}, are uniquely specified for all \( (c',c) \in P_{m,n} \)
  based on the following observation. Pick any \( c' \in \mathcal{C}_{m} \), \( m < n \), such that
  \( W_{m}(c') \subseteq W_{n}(c) \), let \( j \) be the unique index \( 1 \le j \le l \) such that
  \( W_{m}(c') \subseteq W_{i(j)}(c_{j}) \). Suppose \( c' \) does not belong to the maximal family
  of subregions of \( W_{n}(c) \), hence \( m < i(j) \). Using the inductive
  assumption~\eqref{item:coherence-of-steps}, we find that for any \( x \in W_{m}(c') \)
  \begin{displaymath}
    \begin{aligned}
      \phi_{n}(x) &= \phi_{i(j)}(x) + \vec{s}_{i(j),n}(c_{j}, c) \\
      &= \phi_{m}(x) + \vec{s}_{m,i(j)}(c',c_{j}) + \vec{s}_{i(j),n}(c_{j}, c).
    \end{aligned}
  \end{displaymath}
  Thus, for \( s_{m,n}(c',c) = s_{m,i(j)}(c',c_{j}) + s_{i(j),n}(c_{j}, c) \)
  item~\eqref{item:coherence-of-steps} holds for all \( m < n \) and all
  \( c' \in \mathcal{C}_{m} \) such that \( W_{m}(c') \subseteq W_{n}(c) \).

  We check that~\eqref{item:far-from-boundary} continues to hold.  Let us assume that this property
  has been verified for regions at levels below \( n \), and by construction we have also
  established
  \begin{equation}
    \label{eq:2}
      t_{n}(c) - s_{i(j),n}(c_{j},c) \ge t_{i(j)}(c_{j}) + 1.
  \end{equation}
  Using the additivity of values \( s_{m,n}(c',c) \) shown above we get
  \begin{displaymath}
    \begin{aligned}
      t_{n}(c) - t_{m}(c') - s_{m,n}(c',c)
      &= t_{n}(c) - s_{i(j),n}(c_{j},c) - s_{m,i(j)}(c',c_{j}) - t_{m}(c') \\
      \textrm{Eq.~\eqref{eq:2}} &\ge t_{i(j)}(c_{j}) + 1 - s_{m,i(j)}(c',c_{j}) - t_{m}(c') \\
      \textrm{inductive assumption} &\ge t_{m}(c') + 2 - t_{m}(c') = 2.
    \end{aligned}
  \end{displaymath}
  Therefore~\eqref{item:far-from-boundary} holds for all \( m < n \) and all \( (c',c) \in P_{m,n}
  \).

  Finally, to guarantee item~\eqref{item:countably-many-cases} note that diffeomorphisms
  \( \phi_{n} \) have been chosen using Lemma~\ref{lem:moving-smooth-disks} based on the shapes of
  regions \( \widetilde{W}_{n}(c) \), as well as shapes of subregions
  \( \widetilde{W}_{i(j)}(c_{j}) \), and their locations inside \( \widetilde{W}_{n}(c) \) specified
  by values \( \rho(c,c_{j}) \). By Theorem~\ref{thm:ball-shaped-toast}, the union
  \( \bigcup_{k} \mathcal{C}_{k} \) is on a rational grid, so all the values \( \rho(c,c_{j}) \) are
  rational. Also, by~\ref{thm:ball-shaped-toast}\eqref{item:regions-have-countably-many-shapes},
  there are only countably many possible shapes for regions \( W_{i(j)}(c_{j}) \). We may therefore
  choose the same diffeomorphism \( \phi_{n}|_{W_{n}(c)} \) whenever the inputs to
  Lemma~\ref{lem:moving-smooth-disks} are the same, which guarantees fulfillment of
  item~\eqref{item:countably-many-cases}.
\end{proof}

The item~\ref{lem:coherent-diffs}\eqref{item:far-from-boundary} above guarantees that the image of
\( W_{m}(c') \) under \( \phi_{n} \) is at least \( 1 \) unit of distance away from the boundary of
\( B(t_{n}(c)) \) whenever \( W_{m}(c') \subseteq W_{n}(c) \). The following two lemmas show that we
can find such \( n \) and \( c \in \mathcal{C}_{n} \) for which the set \( \phi_{n}(W_{m}(c')) \) is
as far from the boundary of \( B(t_{n}(c)) \) as we desire.

\begin{lemma}
  \label{lem:image-distance-to-boundary}
  Let \( m_{1} < m_{2} < m_{3} < \cdots \) be an increasing sequence of integers and
  \( c_{j} \in \mathcal{C}_{m_{j}} \) be elements such that
  \( W_{m_{j}}(c_{j}) \subseteq W_{m_{j+1}}(c_{j+1}) \) (such a sequence is produced by
  Lemma~\ref{lem:toast-properties}).  For all \( j \ge 2 \) one has
  \[ t_{m_{j}}(c_{j}) - \sum_{1 \le k <j} s_{m_{k},m_{k+1}}(c_{k},c_{k+1}) - t_{m_{1}}(c_{1}) \ge
    j-1. \]
\end{lemma}

\begin{proof}
  The argument is a simple induction coupled with
  item~\ref{lem:coherent-diffs}\eqref{item:far-from-boundary}, which, in particular, gives the base
  \[ t_{m_{2}}(c_{2}) - s_{m_{1},m_{2}}(c_{1},c_{2}) - t_{m_{1}}(c_{1}) \ge 1. \]
  Suppose the statement has been established for \( j-1 \):
  \begin{equation}
    \label{eq:3}
    t_{m_{j-1}}(c_{j-1}) - \mkern-16mu\sum_{1 \le k <j-1}\mkern-16mu s_{m_{k},m_{k+1}}(c_{k},c_{k+1}) - t_{m_{1}}(c_{1}) \ge j-2.
  \end{equation}
  By item~\ref{lem:coherent-diffs}\eqref{item:far-from-boundary} we have
  \( t_{m_{j}}(c_{j}) - s_{m_{j-1},m_{j}}(c_{j-1},c_{j}) \ge 1 + t_{m_{j-1}}(c_{j-1}) \), and
  therefore
  \begin{displaymath}
    \begin{aligned}
      t_{m_{j}}(c_{j}) &- \mkern-10mu\sum_{1 \le k <j}\mkern-10mu s_{m_{k},m_{k+1}}(c_{k},c_{k+1}) - t_{m_{1}}(c_{1}) \\
      &= t_{m_{j}}(c_{j}) - s_{m_{j-1}, m_{j}}(c_{j-1},c_{j}) - \mkern-16mu\sum_{1 \le k <j-1}\mkern-16mu
      s_{m_{k},m_{k+1}}(c_{k},c_{k+1}) - t_{m_{1}}(c_{1}) \\
      &\ge 1 + t_{m_{j-1}}(c_{j-1}) - \mkern-16mu\sum_{1 \le k <j-1}\mkern-16mu
      s_{m_{k},m_{k+1}}(c_{k},c_{k+1}) - t_{m_{1}}(c_{1}) \\
      \textrm{Eq.~\eqref{eq:3}} &\ge 1 + j-2 = j-1,
    \end{aligned}
  \end{displaymath}
  which yields the step of induction.
\end{proof}

\begin{lemma}
  \label{lem:image-far-from-boundary}
  For any \( x \in \Omega \) and any \( R \in \mathbb{R}^{\ge 0} \) there exist \( n \) and
  \( c \in \mathcal{C}_{n} \) such that \( x \in W_{n}(c) \) and 
  \[ ||\phi_{n}(x)|| + R \le t_{n}(c). \]
\end{lemma}
  
\begin{proof}
  In view of~\ref{thm:ball-shaped-toast}\eqref{item:regions-exhaust-space} there is some
  \( m_{1} \) and \( c_{1} \in \mathcal{C}_{m_{1}} \) such that \( x \in W_{m_{1}}(c_{1}) \).  By
  Lemma~\ref{lem:toast-properties} there exist levels \( m_{1} < m_{2} < \cdots \) and points
  \( c_{j} \in \mathcal{C}_{m_{j}} \) such that \( W_{m_{j}}(c_{j}) \subseteq W_{m_{j+1}}(c_{j+1}) \).
  For each \( j \) we have in view of~\ref{lem:coherent-diffs}\eqref{item:coherence-of-steps}
  \begin{displaymath}
    \begin{aligned}
      t_{m_{j}}(c_{j}) - ||\phi_{m_{j}}(x)|| &= t_{m_{j}}(c_{j}) -
      ||\phi_{m_{j-1}}(x) + \vec{s}_{m_{j-1},m_{j}}(c_{j-1},c_{j})|| \\
      &= \cdots\quad \textrm{\( j-2 \) further applications
        of~\ref{lem:coherent-diffs}\eqref{item:coherence-of-steps}}\\
      &= t_{m_{j}}(c_{j}) - \Bigl|\Bigl|\phi_{m_{1}}(x) + \mkern-10mu\sum_{1\le k
        < j}\mkern-7mu\vec{s}_{m_{k},m_{k+1}}(c_{k},c_{k+1}) \Bigr|\Bigr| \\
      &\ge t_{m_{j}}(c_{j}) - t_{m_{1}}(c_{1}) -
      \mkern-10mu\sum_{1\le k < j}\mkern-7mus_{m_{k},m_{k+1}}(c_{k},c_{k+1})\\
      \textrm{Lemma~\ref{lem:image-distance-to-boundary}} &\ge j-1.
    \end{aligned}
  \end{displaymath}
  Therefore \( n = m_{j} \) and \( c = c_{j} \) satisfy the conclusions of the lemma as long as
  \( j-1 \ge R\).
\end{proof}

We now define a new flow \( \mathfrak{F}' \) on the same phase space \( \Omega \). Notation
\( x \oplus \vec{r} \) will be used to distinguish the action of \( \mathfrak{F}' \) from the action
of the original flow \( \mathfrak{F} \). For \( x \in \Omega \) and \( \vec{r} \in \mathbb{R}^{d} \)
let \( c \in \mathcal{C}_{n} \) be such that \( x \in W_{n}(c) \) and
\( \phi_{n}(x) + \vec{r} \in B(t_{n}(c)) \) (such \( n \) and \( c \) exist by
Lemma~\ref{lem:image-far-from-boundary}). The \( \mathfrak{F}' \) action of \( \vec{r} \) upon
\( x \) is defined by
\[ x \oplus \vec{r} = (\phi_{n}|_{W_{n}(c)})^{-1}(\phi_{n}(x) + \vec{r}), \]
or, equivalently, \( x \oplus \vec{r} \) is the element of \( W_{n}(c) \) for which
\( \phi_{n}(x \oplus \vec{r}) = \phi_{n}(x) + \vec{r} \). The geometric interpretation of the action
is as follows. We use the diffeomorphism \( \phi_{n} \) to identify \( W_{n}(c) \) with
\( B(t_{n}(c)) \). One acts upon \( \phi_{n}(x) \in B(t_{n}(c)) \subseteq \mathbb{R}^{d} \) by
translation. Assuming the image lies within the same ball \( B(t_{n}(c)) \), we can pull it back to
an element of \( W_{n}(c) \), which is what \( x \oplus \vec{r} \) is defined to be. As we argue
below, this definition does not depend on the choice of \( n \) and \( c \) due to the coherence of
diffeomorphisms \( \phi_{n} \) provided by
item~\ref{lem:coherent-diffs}\eqref{item:coherence-of-steps}. Having this simple picture of the
action in their mind will make it easy for the reader to follow the somewhat tedious but elementary
computations that constitute a large portion of the remainder of this section.

\begin{lemma}
  \label{lem:oplus-is-well-defined}
  The definition of \( x \oplus \vec{r} \) does not depend on the choice of \( n \) and \( c \).
\end{lemma}

\begin{proof}
  Let \( c' \in \mathcal{C}_{m} \) be another element that can be used in the definition of \( x
  \oplus \vec{r} \), i.e., \( x \in W_{m}(c') \) and
  \( \phi_{m}(x) + \vec{r} \in B(t_{m}(c')) \).
  Item~\ref{thm:ball-shaped-toast}\eqref{item:regions-are-disjoint} implies \( m \ne n \),
  and~\ref{thm:ball-shaped-toast}\eqref{item:regions-are-coherent} guarantees that either
  \( W_{m}(c') \subseteq W_{n}(c) \) or \( W_{n}(c) \subseteq W_{m}(c') \).  Since roles of \( m \)
  and \( n \) are symmetric, we may assume without loss of generality that the former is the
  case. Consider the chain of equalities
  \begin{displaymath}
    \begin{aligned}
      \phi_{n}(x) + \vec{r} &= \phi_{m}(x) + \vec{s}_{m,n}(c',c) + \vec{r} \\
      &= \phi_{m}\Bigl((\phi_{m}|_{W_{m}(c')})^{-1}\bigl(\phi_{m}(x) +
      \vec{r}\bigr)\Bigr) +  \vec{s}_{m,n}(c',c) \\
      \textrm{item~\ref{lem:coherent-diffs}\eqref{item:coherence-of-steps}} 
      &= \phi_{n}\Bigl((\phi_{m}|_{W_{m}(c')})^{-1}\bigl(\phi_{m}(x) +
      \vec{r}\bigr)\Bigr). \\
    \end{aligned}
  \end{displaymath}
  Applying \( (\phi_{n}|_{W_{n}(c)})^{-1} \) to the first and the last expressions above yields
  \[ (\phi_{n}|_{W_{n}(c)})^{-1}(\phi_{n}(x) + \vec{r}) = (\phi_{m}|_{W_{m}(c')})^{-1}(\phi_{m}(x) +
    \vec{r}), \]
  which finishes the proof of the lemma.
\end{proof}

Having established that \( x \oplus \vec{r} \) is well-defined, we can now verify it to be a free Borel
flow.

\begin{lemma}
  \label{lem:oplus-is-a-flow}
  The map \( \Omega \times \mathbb{R}^{d} \ni (x, \vec{r}) \mapsto x \oplus \vec{r} \in \Omega \)
  defines a free Borel action of \( \mathbb{R}^{d} \) on \( \Omega \). This flow is smoothly
  equivalent to \( \mathfrak{F} \). Moreover, the identity map \( \mathrm{id} : \Omega \to \Omega \)
  is a smooth equivalence between \( \mathfrak{F} \) and \( \mathfrak{F}' \).
\end{lemma}

\begin{proof}
  Pick \( \vec{r}_{1}, \vec{r}_{2} \in \mathbb{R}^{d} \), and use
  Lemma~\ref{lem:image-far-from-boundary} to choose \( n \), \( c \in \mathcal{C}_{n} \), such
  that
  \[ ||\phi_{n}(x)|| + ||\vec{r}_{1}|| + ||\vec{r}_{2}|| \le t_{n}(c). \]
  This inequality guarantees that \( (\phi_{n}|_{W_{n}(c)})^{-1} \) is defined in the following terms:
  \begin{displaymath}
    \begin{aligned}
      (x \oplus \vec{r}_{1}) \oplus \vec{r}_{2}
      &=  (\phi_{n}|_{W_{n}(c)})^{-1}\bigl(\phi_{n}(x \oplus \vec{r}_{1}) + \vec{r}_{2}\bigr) \\
      &=  (\phi_{n}|_{W_{n}(c)})^{-1}\bigl(\phi_{n}(x)+ \vec{r}_{1} + \vec{r}_{2}\bigr)
      = x \oplus (\vec{r}_{1} + \vec{r}_{2}). \\
    \end{aligned}
  \end{displaymath}
  Coupled with the straightforward \( x \oplus \vec{0} = x \), these computations show that
  \( \oplus \) defines a flow on \( \Omega \). This flow is free, because the maps
  \( \phi_{n}|_{W_{n}(c)} \) are injective.

  It is easy to see that orbits of \( \mathfrak{F}' \) coincide with those of \( \mathfrak{F} \). The
  inclusion \( E_{\mathfrak{F}'} \subseteq E_{\mathfrak{F}} \) is guaranteed by the condition
  \( W_{n} \subseteq E_{\mathfrak{F}} \). For the inverse direction, let \( x, y \in \Omega \) be
  such that \( x E_{\mathfrak{F}} y \).
  By~\ref{thm:ball-shaped-toast}\eqref{item:regions-exhaust-space}, there are \( n \) and
  \( c \in \mathcal{C}_{n} \) such that \( x, y \in W_{n}(c) \). One has
  \begin{displaymath}
    \begin{aligned}
      x \oplus \bigl(\phi_{n}(y) - \phi_{n}(x)\bigr)
      &= (\phi_{n}|_{W_{n}(c)})^{-1}\bigl(\phi_{n}(x) + \phi_{n}(y) - \phi_{n}(x)\bigr) \\
      &= (\phi_{n}|_{W_{n}(c)})^{-1}(\phi_{n}(y)) = y,
    \end{aligned}
  \end{displaymath}
  and thus \( E_{\mathfrak{F}} = E_{\mathfrak{F}'} \).

  The flow \( \mathfrak{F}' \) is Borel. To justify this set
  \( \alpha : \Omega \times \mathbb{R}^{d} \to \mathbb{R}^{d} \) to be defined by
  \( x \oplus \vec{r} = x + \alpha(x, \vec{r}) \); it suffices to show that \( \alpha \) is Borel.
  In general, we would have to verify that the value
  \( (\phi_{n}|_{W_{n}(c)})^{-1}\bigl(\phi_{n}(x) + \vec{r}\,\bigr) \) depends in a Borel way on
  \( x \in \Omega \) and \( \vec{r} \in \mathbb{R}^{d} \), which requires going into the details of
  the way maps \( \phi_{n}|_{W_{n}(c)} \) are constructed based on \( n \) and
  \( c \in \mathcal{C}_{n} \). However, we ensured in
  item~\ref{lem:coherent-diffs}\eqref{item:countably-many-cases} that there is a countable Borel
  partition of each cross section \( \mathcal{C}_{n} = \bigsqcup_{k}\mathcal{C}_{n,k} \) such that
  \( \widetilde{W}_{n}(c_{1}) = \widetilde{W}_{n}(c_{2}) \), \( t_{n}(c_{1}) = t_{n}(c_{2}) \), and
  \( \phi_{n}|_{W_{n}(c_{1})}(c_{1} + \vec{r}) = \phi_{n}|_{W_{n}(c_{2})}(c_{2} + \vec{r})\) for all
  \( c_{1},c_{2} \in \mathcal{C}_{n,k}\) and all \( \vec{r} \in \widetilde{W}(c_{1}) \). Let
  \( \widetilde{W}_{n,k} \) denote the common shape of \( \widetilde{W}_{n}(c) \), and similarly let
  \( t_{n,k} \) be the radius \( t_{n}(c) \) common for all \( c \in \mathcal{C}_{n,k} \). Likewise,
  maps \( \phi_{n}|_{W_{n}(c)}(c + \cdot ) \) produce the same diffeomorphism
  \( \varphi_{n,k} : \widetilde{W}_{n,k} \to B(t_{n,k}) \) regardless of the choice of
  \( c \in \mathcal{C}_{n,k} \).

  Recall that \( \pi_{n} : V_{n} \to \mathcal{C}_{n} \) is the map that associates the distinguished
  point \( c \) to every \( x \) in \( W_{n}(c) \). Let
  \[ X_{n,k} = \{ (x, \vec{r}) \in \Omega \times \mathbb{R}^{d} : \pi_{n}(x) \in \mathcal{C}_{n,k},\
    x \in V_{n}, \textrm{ and } \varphi_{n,k}(\rho(\pi_{n}(x), x)) + \vec{r} \le t_{n,k}\} \]
  denote the set of those pairs \( (x, \vec{r}) \) for which \( x \oplus \vec{r} \) can be defined
  using \( W_{n}(c) \) for some \( c \in \mathcal{C}_{n,k} \). To conclude that \( \alpha \) is
  Borel, we observe that for \( (x, \vec{r}) \in X_{n,k} \) its value equals
  \begin{equation}
    \label{eq:8}
    \alpha(x, \vec{r}) = \rho(x, \pi_{n}(x)) +
    \varphi_{n,k}^{-1}\bigl(\varphi_{n,k}(\rho(\pi_{n}(x),x))
    + \vec{r} \bigr),
  \end{equation}
  which is a composition of Borel functions. It is at this point that our efforts in ensuring
  stability of the construction
  in~\ref{thm:ball-shaped-toast}\eqref{item:regions-have-countably-many-shapes}
  and~\ref{lem:coherent-diffs}\eqref{item:countably-many-cases} yield their fruits. We have only one
  diffeomorphism \( \varphi_{n,k} \) in the definition of \( \alpha \) which applies to all
  arguments \( (x, \vec{r}) \in X_{n,k} \). Thus \( \alpha|_{X_{n,k}} \) is Borel regardless of how
  this diffeomorphism was chosen. Since \( \bigcup_{n,k} X_{n,k} = \Omega \times \mathbb{R}^{d} \),
  we may conclude that \( \alpha \) is Borel on all of \( \Omega \times \mathbb{R}^{d} \).

  We have already established that \( \mathrm{id} : \Omega \to \Omega \) is an orbit equivalence,
  and it remains to verify that it is smooth, which amounts to showing that
  \( \alpha(x, \cdot) : \mathbb{R}^{d} \to \mathbb{R}^{d} \) is a diffeomorphism for each
  \( x \in \Omega \). It is a bijection, since \( \mathfrak{F},\ \mathfrak{F}' \) are free and
  \( E_{\mathfrak{F}} = E_{\mathfrak{F}'} \), and it is smooth and orientation preserving, since
  according to Eq.~\eqref{eq:8} for any fixed \( x \in \Omega \) the function \( \alpha(x, \cdot) \)
  is a composition of \( \varphi_{n,k}^{-1} \) with translation maps.
\end{proof}

The goal of this section is to show that every free flow is smoothly equivalent to a special one. So
far starting with a free flow \( \mathfrak{F} \) we have constructed a smoothly equivalent flow
\( \mathfrak{F}' \), and it remains to verify that the latter is special. We need a subset
\( \Omega_{1} \subseteq \Omega \) invariant under the shifts
\( \Omega_{1} \oplus s \times \vec{0}^{d-1} = \Omega_{1} \) for all \( s \in \mathbb{R} \). Set
\( L_{n} \subseteq W_{n} \) to correspond to the preimages of the \( x \)-axis inside
\( B(t_{n}(c)) \) under the diffeomorphisms \( \phi_{n} \):
\[ L_{n} = \bigl\{ (c, x) \in W_{n} : \phi_{n}(x) \in \bigl[-t_{n}(c), t_{n}(c)\bigr] \times
  \vec{0}^{d-1}\bigr\}. \]
Sets \( L_{n}(c) \) represent line segments inside regions \( W_{n}(c) \)
(see~Figure~\ref{fig:inductive-step}). Set \( L = \bigcup_{n}L_{n} \), and let \( \Omega_{1} \) be
the projection of \( L \) onto the second coordinate:
\[ \Omega_{1} = \Bigl\{x : (c, x) \in L \textrm{ for some } c \in \bigcup_{n} \mathcal{C}_{n} \Bigr\}. \]

\begin{lemma}
  \label{lem:lset-invariant-first-coordinate}
  The set \( \Omega_{1} \subseteq \Omega \) is Borel and
  \( \Omega_{1} \oplus s \times \vec{0}^{d-1} = \Omega_{1} \) for all \( s \in \mathbb{R} \).
\end{lemma}

\begin{proof}
  The set \( \Omega_{1} = \bigcup_{n} \proj_{2}(L_{n}) \) is Borel, since projections
  \( \proj_{2} : L_{n} \to \Omega_{1} \) are injective, hence have Borel images
  (see~\cite[Corollary~15.2]{kechris_classical_1995}). To check shift invariance, pick
  \( x \in \Omega_{1} \) and \( s \in \mathbb{R} \); let \( \vec{s} \) denote the vector
  \( s \times \vec{0}^{d-1} \). There has to exist some \( m \) and \( c' \in \mathcal{C}_{m} \)
  such that \( (c',x) \in L_{m} \). Pick \( n > m \) and \( c \in \mathcal{C}_{n} \) to satisfy
  \( x, x \oplus \vec{s} \in W_{n}(c) \), which exist
  by~\ref{thm:ball-shaped-toast}\eqref{item:regions-exhaust-space}. Note that \( (c,x) \in L_{n} \)
  as according to~\ref{lem:coherent-diffs}\eqref{item:coherence-of-steps}
  \( \phi_{n}(x) = \phi_{m}(x) + \vec{s}_{m,n}(c',c), \) and therefore
  \( \phi_{m}(x) \in [-t_{m}(c'), t_{m}(c')] \times \vec{0}^{d-1} \) implies
  \begin{displaymath}
    \begin{aligned}
      \phi_{n}(x) &\in \bigl[-t_{m}(c') - s_{m,n}(c',c),\ t_{m}(c')+s_{m,n}(c',c)\bigr]
      \times
      \vec{0}^{d-1} \\
      \textrm{item~\ref{lem:coherent-diffs}\eqref{item:far-from-boundary}}&\subseteq [-t_{n}(c),
      t_{n}(c)] \times \vec{0}^{d-1}.
    \end{aligned}
  \end{displaymath}
  Since \( x \oplus \vec{s} \in W_{n}(c) \), we have
  \( \phi_{n}(x) + \vec{s} \in [-t_{n}(c), t_{n}(c)] \times \vec{0}^{d-1} \), thus
  \( x \oplus \vec{s} \in \Omega_{1} \) as claimed.
\end{proof}

The following lemma will be helpful in establishing that \( \mathfrak{F}' \) is special.

\begin{lemma}
  \label{lem:projection-F-prime}
  For any \( x \in \Omega \) there exists \( \vec{q} \in \mathbb{R}^{d-1} \) such that \( x \oplus 0
  \times \vec{q} \in \Omega_{1} \).
\end{lemma}

\begin{proof}
  Pick some \( x \in \Omega \) and, as usual, let \( n \), \( c \in \mathcal{C}_{n} \) be chosen to
  satisfy \( x \in W_{n}(c) \).  Set \( \vec{q} \in \mathbb{R}^{d-1} \) to be the negative of the
  projection of \( \phi_{n}(x) \) onto the last \( (d-1) \)-many coordinates:
  \( \vec{q} = -\proj_{[2,d]}(\phi_{n}(x)) \).  By the definition of the
  action,
  \begin{displaymath}
    \begin{aligned}
      x \oplus 0 \times \vec{q} &= (\phi_{n}|_{W_{n}(c)})^{-1}(\phi_{n}(x) + 0 \times \vec{q}\,) \\
      &= (\phi_{n}|_{W_{n}(c)})^{-1}\bigl(\phi_{n}(x) - 0 \times
      \proj_{[2,d]}(\phi_{n}(x))\bigr)  \\
      &= (\phi_{n}|_{W_{n}(c)})^{-1}\bigl(\proj_{1}(\phi_{n}(x)) \times \vec{0}^{d-1}\bigr) \in
      L_{n}(c) \subseteq \Omega_{1},
    \end{aligned}
  \end{displaymath}
  whence \( \vec{q} \) satisfies the conclusion of the lemma.
\end{proof}

We have established that \( \Omega_{1} \) is invariant under the \( \mathbb{R} \)-flow corresponding
to the actions by vectors \( s \times \vec{0} \), and we may therefore naturally define a special
flow \( \overline{\mathfrak{F}}' \) on \( \Omega_{1} \times \mathbb{R}^{d-1} \) by
\[ (x, \vec{q}\,) \boxplus (r_{1}, r_{2}, \ldots, r_{d}) = \bigl(x \oplus r_{1} \times \vec{0}^{d-1},\
  \vec{q} + (r_{2}, \ldots, r_{d})\bigr).\]
Note that \( \boxplus \) is used for the action to distinguish it from both the actions given by
\( \mathfrak{F} \) and \( \mathfrak{F}' \).  We are going to verify that
\( \overline{\mathfrak{F}}' \) is isomorphic to \( \mathfrak{F}' \), and to this end we define two
Borel maps \( \mu : \Omega \to \Omega_{1} \) and \( \nu : \Omega \to \mathbb{R}^{d-1} \) such that
\( \Omega \ni x \mapsto (\mu(x), \nu(x)) \in \Omega_{1} \times \mathbb{R}^{d-1} \) will be the
desired isomorphism.  For \( x \in \Omega \), \( n \) and \( c \in \mathcal{C}_{n} \),
\( x \in W_{n}(c) \), set

\begin{displaymath}
  \begin{aligned}
   \mu(x) &= (\phi_{n}|_{W_{n}(c)})^{-1}\bigl(\proj_{1}(\phi_{n}(x)) \times \vec{0}^{d-1}\bigr) \in
   \Omega_{1} \\
   \nu(x) &= \proj_{[2,d]}(\phi_{n}(x)) \in \mathbb{R}^{d-1}.
  \end{aligned}
\end{displaymath}

\begin{lemma}
  \label{lem:psi-is-well-defined}
  Maps \( \mu \) and \( \nu \) are well-defined in the sense that their values do not depend on the
  choice of \( n \) and \( c \in \mathcal{C}_{n} \).
\end{lemma}

\begin{proof}
  Let \( m < n \) and \( c' \in \mathcal{C}_{m} \) be other elements such that \( x \in W_{m}(c) \).
  Once again, \( \phi_{n}(x) = \phi_{m}(x) + \vec{s}_{m,n}(c',c) \), by
  item~\ref{lem:coherent-diffs}\eqref{item:coherence-of-steps}. In particular, the projections of
  \( \phi_{n}(x) \) and \( \phi_{m}(x) \) onto the last \( (d-1) \)-many coordinates are equal,
  because \( \vec{s}_{m,n}(c',c) = s_{m,n}(c',c) \times \vec{0}^{d-1} \) This shows that
  \( \nu(x) \) is well-defined.

  For \( s = \proj_{1}(\phi_{n}(x)) \), \( s' = \proj_{1}(\phi_{m}(x)) \), and
  \( \vec{s} = s \times \vec{0}^{d-1} \), \( \vec{s}\,' = s' \times \vec{0}^{d-1}\), we have
  \begin{displaymath}
    \begin{aligned}
      \vec{s} &= \vec{s}\,' + \vec{s}_{m,n}(c',c) \\
      &= \phi_{m}\bigl((\phi_{m}|_{W_{m}(c')})^{-1}(\vec{s}\,')\bigr) + \vec{s}_{m,n}(c',c) \\
      \textrm{item~\ref{lem:coherent-diffs}\eqref{item:coherence-of-steps}}&=
      \phi_{n}\bigl((\phi_{m}|_{W_{m}(c')})^{-1}(\vec{s}\,')\bigr).
    \end{aligned}
  \end{displaymath}
  Applying \( (\phi_{n}|_{W_{n}(c)})^{-1} \) to both sides yields
  \begin{displaymath}
    \begin{aligned}
      (\phi_{n}|_{W_{n}(c)})^{-1}(\vec{s}) = (\phi_{m}|_{W_{m}(c')})^{-1}(\vec{s}\,'),
    \end{aligned}
  \end{displaymath}
  which shows that the value \( \mu(x) \in \Omega_{1} \) does not depend on the choice of \( n \)
  and \( c \in \mathcal{C}_{n} \).
\end{proof}

\begin{lemma}
  \label{lem:psi-is-bijective}
  The map \( \Omega \ni x \mapsto (\mu(x), \nu(x)) \in \Omega_{1} \times \mathbb{R}^{d-1} \) is a
  bijection.
\end{lemma}

\begin{proof}
  For injectivity, let \( x, y \in \Omega \) be distinct; recall that the orbit equivalence
  relations of the flows \( \mathfrak{F} \) and \( \mathfrak{F}' \) coincide by
  Lemma~\ref{lem:oplus-is-a-flow}, and we denote it by \( E \). Note that \( x E \mu(x) \) and
  \( y E \mu(y) \), so if \( \neg x E y \), then \( \mu(x) \ne \mu(y) \). Thus we need to consider
  the case \( x E y \) and by~\ref{thm:ball-shaped-toast}\eqref{item:regions-exhaust-space} there is
  some \( n \) and \( c \in \mathcal{C}_{n} \) such that \( x, y \in W_{n}(c) \). Since
  \( \phi_{n}|_{W_{n}(c)} : W_{n}(c) \to B(t_{n}(c)) \) is injective,
  \( \phi_{n}(x) \ne \phi_{n}(y) \), and thus either
  \begin{displaymath}
    \begin{aligned}
      \proj_{1}(\phi_{n}(x)) &\ne \proj_{1}(\phi_{n}(y)) \qquad \textrm{or}\\
      \proj_{[2,d]}(\phi_{n}(x)) &\ne \proj_{[2,d]}(\phi_{n}(y)),
    \end{aligned}
  \end{displaymath}
  and hence either \( \mu(x) \ne \mu(y) \) or \( \nu(x) \ne \nu(y) \).

  For surjectivity, pick \( (x, \vec{q}\,) \in \Omega_{1} \times \mathbb{R}^{d-1} \), and
  \( n \), \( c \in \mathcal{C}_{n} \), with \( x \in W_{n}(c) \) and
  \( t_{n}(c) - ||\phi_{n}(x)|| \ge ||\vec{q}|| \), which exist by
  Lemma~\ref{lem:image-far-from-boundary}.  Note that by the coherence property of \( L_{n} \)
  established in the proof of Lemma~\ref{lem:lset-invariant-first-coordinate}, we have
  \( (c,x) \in L_{n} \), which is equivalent to \( \phi_{n}(x) \in \mathbb{R} \times \vec{0}^{d-1} \).
  These conditions ensure that \( (\phi_{n}|_{W_{n}(c)})^{-1} \) is defined on
  \( \phi_{n}(x) + 0 \times \vec{q} \) and
  \[ (\mu,\nu)\bigl((\phi_{n}|_{W_{n}(c)})^{-1}\bigl(\phi_{n}(x) + 0 \times \vec{q}\,\bigr)\bigr) =
    (x, \vec{q}\,),\]
  witnessing surjectivity.
\end{proof}

At last, we can check that the flows \( \mathfrak{F}' \) and \( \overline{\mathfrak{F}}' \) are
isomorphic.

\begin{lemma}
  \label{lem:psi-isomorphism-oplus}
  The map \( \Omega \ni x \mapsto (\mu(x), \nu(x)) \in \Omega_{1} \times \mathbb{R}^{d-1} \) is an
  isomorphism of flows \( \mathfrak{F}' \) and \( \overline{\mathfrak{F}}' \).
\end{lemma}

\begin{proof}
  Once Lemma~\ref{lem:psi-is-bijective} is available, our remaining goal is to show that
  \( (\mu, \nu)(x \oplus \vec{r}\,) = (\mu(x), \nu(x)) \boxplus \vec{r} \) for all
  \( x \in \Omega \) and all \( \vec{r} \in \mathbb{R}^{d} \). We first verify this for those
  \( \vec{r} \) that satisfy \( \proj_{1}(\vec{r}\,) = 0 \), i.e., \( \vec{r} = 0 \times \vec{q} \),
  for some \( \vec{q} \in \mathbb{R}^{d-1} \). One has
  \begin{displaymath}
    \begin{aligned}
      \nu(x \oplus 0 \times \vec{q}\,) &= \proj_{[2,d]}\bigl(\phi_{n}(x) + 0\times \vec{q}\,\bigr) \\
      & = \proj_{[2,d]}\bigl(\phi_{n}(x) \bigr) + \vec{q} = \nu(x) + \vec{q}, \\
      \mu(x \oplus 0 \times \vec{q}\,) &= (\phi_{n}|_{W_{n}(c)})^{-1}\bigl(\proj_{1}\bigl(\phi_{n}(x) +
      0\times \vec{q}\,\bigr)\times \vec{0}^{d-1}\bigr) \\
      &= (\phi_{n}|_{W_{n}(c)})^{-1}\bigl(\proj_{1}(\phi_{n}(x)) \times \vec{0}^{d-1} \bigr) = \mu(x).
    \end{aligned}
  \end{displaymath}
  We have thus shown that
  \begin{equation}
    \label{eq:6}
    (\mu,\nu)(x \oplus 0 \times \vec{q}\,) = \bigl((\mu,\nu)(x)\bigr) \boxplus 0 \times \vec{q}
    \quad \textrm{for all \( \vec{q} \in \mathbb{R}^{d-1} \) and all \( x \in \Omega \).}
  \end{equation}
  Note also that \( (\mu,\nu)(y) = (y,\vec{0}) \) for all \( y \in \Omega_{1} \subseteq \Omega \),
  and therefore by Lemma~\ref{lem:lset-invariant-first-coordinate} for all \( s \in \mathbb{R} \)
  and all \( y \in \Omega_{1} \)
  \begin{equation}
    \label{eq:7}
    (\mu,\nu)(y \oplus s \times \vec{0}^{d-1}) = (y \oplus s \times \vec{0}^{d-1}, \vec{0}) = (\mu,\nu)(y)
    \boxplus (s \times \vec{0}^{d-1}).
  \end{equation}

  For any \( x \in \Omega \), Lemma~\ref{lem:projection-F-prime} gives an element
  \( \vec{q}_{0} \in \mathbb{R}^{d-1} \) such that
  \( x \oplus 0 \times \vec{q}_{0} \in \Omega_{1} \). Let \( \vec{r} \in \mathbb{R}^{d} \) be
  arbitrary, and write it as \( \vec{r} = s \times \vec{0}^{d-1} + 0 \times \vec{q} \) for some
  \( s \in \mathbb{R}\) and \( \vec{q} \in \mathbb{R}^{d-1} \).  We have
  \begin{displaymath}
    \begin{aligned}
      (\mu,\nu)(x \oplus \vec{r}\,) &= (\mu,\nu)\bigl(x \oplus 0 \times \vec{q}_{0} \oplus (-0 \times
      \vec{q}_{0}) \oplus s \times
      \vec{0}^{d-1} \oplus 0 \times \vec{q}\,\bigr) \\
      &= (\mu,\nu)\bigl(x \oplus 0 \times \vec{q}_{0} \oplus s \times
      \vec{0}^{d-1} \oplus 0 \times (\vec{q}- \vec{q}_{0})\bigr) \\
      \textrm{Eq.~\eqref{eq:6}} &= (\mu,\nu)\bigl(x \oplus 0 \times \vec{q}_{0} \oplus
      s \times \vec{0}^{d-1}\bigr) \boxplus 0 \times (\vec{q} - \vec{q}_{0}) \\
      \textrm{Eq.~\eqref{eq:7}}&= (\mu,\nu)(x \oplus 0 \times \vec{q}_{0}) \boxplus s \times \vec{0}^{d-1}
      \boxplus 0 \times (\vec{q} - \vec{q}_{0}) \\
      \textrm{Eq.~\eqref{eq:6}}&=(\mu,\nu)(x) \boxplus 0 \times \vec{q}_{0} \boxplus s \times
      \vec{0}^{d-1}
      \boxplus 0 \times (\vec{q} - \vec{q}_{0}) \\
      &=(\mu,\nu)(x) \boxplus (s \times \vec{q}\,) = (\mu,\nu)(x) \boxplus \vec{r}.
    \end{aligned}
  \end{displaymath}
  Thus \( (\mu,\nu) \) is an isomorphism between flows \( \mathfrak{F}' \) and
  \( \overline{\mathfrak{F}}' \).
\end{proof}
  
The following theorem summarizes the analysis that has been conducted in
Lemmas~\ref{lem:coherent-diffs} through~\ref{lem:psi-isomorphism-oplus}.

\begin{theorem}
  \label{thm:special-flow-equivalence}
  Every free Borel \( \mathbb{R}^{d} \)-flow, \( d \ge 2 \), is smoothly equivalent to a special
  flow.
\end{theorem}

\section{Smooth Equivalence of Flows}
\label{sec:smooth-equiv-flows}

We are finally ready for the proof of the main result of this article\textemdash{}smooth equivalence of all
non-tame Borel \( \mathbb{R}^{d} \)-flows.  For this we just need to combine
Theorem~\ref{thm:special-flow-equivalence} with the result of Miller--Rosendal on one-dimensional
flows.

\begin{theorem}
  \label{thm:time-change-equivalence}
  All non-tame free Borel \( \mathbb{R}^{d} \)-flows, \( d \ge 2 \), are smoothly equivalent.
\end{theorem}

\begin{proof}
  Let \( \mathfrak{F}_{1} \) and \( \mathfrak{F}_{2} \) be non-tame free Borel
  \( \mathbb{R}^{d} \)-flows.  By Theorem~\ref{thm:special-flow-equivalence}, each of them is
  smoothly equivalent to a special flow
  \( \mathbb{R}^{d} \acts \Omega_{i} \times \mathbb{R}^{d-1} \), \( i = 1,2 \).  Note that neither
  of the \( \mathbb{R} \)-flows \( \mathbb{R} \acts \Omega_{i} \) can be tame, for otherwise the
  corresponding \( \mathbb{R}^{d} \)-flow would also be tame.  By the Miller--Rosendal result
  (see~Theorem~\ref{thm:non-tame-flows-smooth-equivlaence}), there is a smooth equivalence
  \( \bar{\xi} : \Omega_{1} \to \Omega_{2} \) between the \( \mathbb{R} \)-flows.  Let
  \( \alpha_{\bar{\xi}} : \Omega_{1} \times \mathbb{R} \to \mathbb{R} \) be the corresponding family
  of diffeomorphisms defined by \( \alpha_{\bar{\xi}}(x, s) = \rho(\bar{\xi}(x), \bar{\xi}(x+s)) \).
  The map \( \xi : \Omega_{1} \times \mathbb{R}^{d-1} \to \Omega_{2} \times \mathbb{R}^{d-1} \)
  defined by \( \xi(x, \vec{q}\,) = (\bar{\xi}(x), \vec{q}\,) \) is a smooth equivalence between the
  special flows, because
  \[ \alpha_{\xi}\bigl((x, \vec{q}\,), (r_{1}, \ldots, r_{d})\bigr) = (\alpha_{\bar{\xi}}(x, r_{1}),
    r_{2}, \ldots, r_{d}) \in \mathbb{R}^{d}\]
  is a \( C^{\infty} \)-smooth orientation preserving diffeomorphism for all \( x \in \Omega_{1} \)
  and all \( \vec{q} \in \mathbb{R}^{d-1} \).  Flows \( \mathfrak{F}_{1} \) and
  \( \mathfrak{F}_{2} \) are smoothly equivalent by transitivity of the smooth equivalence relation.
\end{proof}

Our approach to the construction of smooth equivalence between multidimensional flows is different
from those taken in the ergodic theoretical antecedents.  Of particular interest is the technique
used in~\cite{MR1113569}.  Their strategy is to start with an orbit equivalence between cross
sections of flows\footnote{The conditions of when such an orbit equivalence exists are well
  understood and follow from Dye's Theorem~\cite{MR131516,MR158048} and Dougherty--Jackson--Kechris
  classification of the hyperfinite equivalence relations~\cite{MR1149121}.}, and extend such a map
to a smooth equivalence by a back-and-forth argument.  It would be interesting to establish the
possibility of such an extension in the descriptive set-theoretical context.  In~\cite{MR3984276},
such an extension was constructed up to a compressible set.  It was also proven therein that if a
complete extension is always possible, then it will imply smooth equivalence of all flows
(Theorem~\ref{thm:time-change-equivalence} above).

\begin{conjecture}
  \label{conj:extending-orbit-equivalences}
  Let \( \mathbb{R}^{d} \acts \Omega_{1} \) and \( \mathbb{R}^{d} \acts \Omega_{2} \) be free Borel
  flows, \( d \ge 2 \), let \( \mathcal{C}_{i} \subseteq \Omega_{i} \) be cocompact cross sections,
  and let \( \zeta : \mathcal{C}_{1} \to \mathcal{C}_{2} \) be an orbit equivalence (i.e., a Borel
  bijection such that \( c E_{1} c' \iff \zeta(c) E_{2} \zeta(c') \)).  There exists a smooth
  equivalence \( \xi : \Omega_{1} \to \Omega_{2} \) that extends \( \zeta \).
\end{conjecture}

Rudolph's Theorem~\cite[Proposition~1.1]{MR536948} produces a smooth equivalence of
\( \mathbb{R}^{d} \)-flows which is also a Lebesgue orbit equivalence, i.e., a map that preserves
the Lebesgue measure between orbits (see~\cite{MR3681992} for the discussion of the notion of
Lebesgue orbit equivalence). This is a significant strengthening of the smooth equivalence, and it
is therefore interesting to see if such a strengthening is possible in the descriptive
set-theoretical context as well. Any Lebesgue orbit equivalence produces an isomorphism between the
spaces of invariant measures of the flows~\cite[Theorem~4.5]{MR3681992}, so we need to restrict
ourselves to flows that have the same cardinalities of sets of ergodic invariant probability
measures. According to~\cite[Theorem~9.1]{MR3681992}, this would ensure the existence of a Lebesgue
orbit equivalence between the flows and, of course, Theorem~\ref{thm:time-change-equivalence}
guarantees that there is a smooth equivalence as well. The question is whether the two can be
combined.\footnote{We would like to thank the referee for suggesting this question.}

\begin{question}
  Let \( \mathfrak{F}_{1} \) and \( \mathfrak{F}_{2} \) be free non-tame Borel
  \( \mathbb{R}^{d} \)-flows, \( d \ge 2 \), having the same numbers of invariant ergodic
  probability measures. Is there a smooth equivalence between these flows which is also a Lebesgue
  orbit equivalence?
\end{question}

\bibliographystyle{halpha}
\bibliography{refs}

\begin{thebibliography}{ORW82}

\bibitem[BJ07]{MR2322367}
Charles~M. Boykin and Steve Jackson.
\newblock Borel boundedness and the lattice rounding property.
\newblock In {\em Advances in logic}, volume 425 of {\em Contemp. Math.}, pages
  113--126. Amer. Math. Soc., Providence, RI, 2007.

\bibitem[DJK94]{MR1149121}
Randall Dougherty, Steve Jackson, and Alexander~S. Kechris.
\newblock The structure of hyperfinite {B}orel equivalence relations.
\newblock {\em Trans. Amer. Math. Soc.}, 341(1):193--225, 1994.

\bibitem[Dye59]{MR131516}
Henry~A. Dye.
\newblock On groups of measure preserving transformations. {I}.
\newblock {\em Amer. J. Math.}, 81:119--159, 1959.

\bibitem[Dye63]{MR158048}
Henry~A. Dye.
\newblock On groups of measure preserving transformations. {II}.
\newblock {\em Amer. J. Math.}, 85:551--576, 1963.

\bibitem[Fel76]{MR409763}
Jacob Feldman.
\newblock New {$K$}-automorphisms and a problem of {K}akutani.
\newblock {\em Israel J. Math.}, 24(1):16--38, 1976.

\bibitem[Fel91]{MR1113569}
Jacob Feldman.
\newblock Changing orbit equivalences of {${\bf R}^d$} actions, {$d\geq 2$}, to
  be {${C}^\infty$} on orbits.
\newblock {\em Internat. J. Math.}, 2(4):409--427, 1991.

\bibitem[Fel92]{MR1163729}
Jacob Feldman.
\newblock Correction to: ``{C}hanging orbit equivalences of {${\bf R}^d$}
  actions, {$d\geq 2$}, to be {$C^\infty$} on orbits''.
\newblock {\em Internat. J. Math.}, 3(3):349--350, 1992.

\bibitem[GJ15]{MR3359054}
Su~Gao and Steve Jackson.
\newblock Countable abelian group actions and hyperfinite equivalence
  relations.
\newblock {\em Invent. Math.}, 201(1):309--383, 2015.

\bibitem[GJKS]{GaoJacksonDiskShapedRegions}
Su~Gao, Steve Jackson, Edward Krohne, and Brandon Seward.
\newblock Borel combinatorics of abelian group actions.
\newblock \textit{In preparation}.

\bibitem[JKL02]{MR1900547}
Steve Jackson, Alexander~S. Kechris, and Alain Louveau.
\newblock Countable {B}orel equivalence relations.
\newblock {\em J. Math. Log.}, 2(1):1--80, 2002.

\bibitem[Kak43]{MR14222}
Shizuo Kakutani.
\newblock Induced measure preserving transformations.
\newblock {\em Proc. Imp. Acad. Tokyo}, 19:635--641, 1943.

\bibitem[Kat75]{MR0412383}
Anatole~B. Katok.
\newblock Time change, monotone equivalence, and standard dynamical systems.
\newblock {\em Dokl. Akad. Nauk SSSR}, 223(4):789--792, 1975.

\bibitem[Kat77]{MR0442195}
Anatole~B. Katok.
\newblock Monotone equivalence in ergodic theory.
\newblock {\em Izv. Akad. Nauk SSSR Ser. Mat.}, 41(1):104--157, 231, 1977.

\bibitem[Kec95]{kechris_classical_1995}
Alexander~S. Kechris.
\newblock {\em Classical descriptive set theory}, volume 156 of {\em Graduate
  Texts in Mathematics}.
\newblock Springer-Verlag, New York, 1995.

\bibitem[MR10]{MR2578608}
Benjamin~D. Miller and Christian Rosendal.
\newblock Descriptive {K}akutani equivalence.
\newblock {\em J. Eur. Math. Soc. (JEMS)}, 12(1):179--219, 2010.

\bibitem[MU17]{MR3702673}
Andrew~S. Marks and Spencer~T. Unger.
\newblock Borel circle squaring.
\newblock {\em Ann. of Math. (2)}, 186(2):581--605, 2017.

\bibitem[Nak88]{MR955378}
Munetaka Nakamura.
\newblock Time change and orbit equivalence in ergodic theory.
\newblock {\em Hiroshima Math. J.}, 18(2):399--412, 1988.

\bibitem[ORW82]{MR653094}
Donald~S. Ornstein, Daniel Rudolph, and Benjamin Weiss.
\newblock Equivalence of measure preserving transformations.
\newblock {\em Mem. Amer. Math. Soc.}, 37(262):xii+116, 1982.

\bibitem[Rud79]{MR536948}
Daniel Rudolph.
\newblock Smooth orbit equivalence of ergodic {${\bf R}^{d}$} actions, {$d\geq
  2$}.
\newblock {\em Trans. Amer. Math. Soc.}, 253:291--302, 1979.

\bibitem[Slu17]{MR3681992}
Konstantin Slutsky.
\newblock Lebesgue orbit equivalence of multidimensional {B}orel flows: a
  picturebook of tilings.
\newblock {\em Ergodic Theory Dynam. Systems}, 37(6):1966--1996, 2017.

\bibitem[Slu19]{MR3984276}
Konstantin Slutsky.
\newblock On time change equivalence of {B}orel flows.
\newblock {\em Fund. Math.}, 247(1):1--24, 2019.

\end{thebibliography}

\end{document}